\documentclass[a4paper, 11pt]{article} 
\usepackage{amsfonts,amsmath,amssymb,amscd}
\usepackage{latexsym}
\usepackage{graphicx}
\usepackage{color} 

\usepackage{hyperref}

\DeclareMathOperator{\Tr}{Tr}

\newtheorem{thm}{Theorem}[section]
\newtheorem{cor}[thm]{Corollary}
\newtheorem{lem}[thm]{Lemma}

\newtheorem{Proposition}[thm]{Proposition}

\newtheorem{defin}[thm]{Definition}
\newtheorem{rem}[thm]{Remark}
\newtheorem{exa}[thm]{Example}

\numberwithin{equation}{section}

\newcommand{\beqr}{\begin{eqnarray}}
\newcommand{\eeqr}{\end{eqnarray}}

\newcommand{\beq}{\begin{eqnarray*}}
\newcommand{\eeq}{\end{eqnarray*}}

\newcommand{\bq}{\begin{equation}}
\newcommand{\eq}{\end{equation}}

\newenvironment{preuve}[1][]
{\vskip 2mm  {\it \bf Proof#1. }}{$\Box$ \vskip 2mm}

\newcommand{\bpr}{\begin{preuve}}
\newcommand{\epr}{\end{preuve}}

\newcommand{\bh}{\mathcal{H}}
\newcommand{\bi}{\mathcal{I}}
\newcommand{\bj}{\mathcal{J}}

\newcommand{\Nn}{\mathbb{N}}
\newcommand{\R}{\mathbb{R}}
\newcommand{\Rp}{\mathbb{R}^*_+}

\newcommand{\C}{\mathbb{C}}

\newcommand{\si}{\sigma}

\newcommand{\dxp}{d_{|x}p }

\newcommand{\gme}{\Gamma(M,E)}

\newcommand{\jle}{\mathcal J^l(E)}
\newcommand{\bjlh}{\bj^l(E_{| \bh})}
\newcommand{\jme}{\mathcal J^m(E)}
\newcommand{\jleh}{\mathcal J^l(E_{|\mathcal H})}
\newcommand{\jmeh}{\mathcal J^m(E_{|\mathcal H})}

\newcommand{\mcp}{M\setminus{Crit(p)}}

\newcommand{\us}{\underline{U}}
\newcommand{\xs}{(x,s)}
\newcommand{\xd}{\dot{x}}
\newcommand{\sd}{\dot{s}}
\newcommand{\xsd}{(\xd, \sd)}
\newcommand{\nae}{\nabla^E}
\newcommand{\suz}{s^{-1}(0)}
\newcommand{\cg}{\langle}
\newcommand{\cd}{\rangle}
\newcommand{\el}{\tilde e_\lambda}
\newcommand{\izn}{i\in \{0, \cdots, n-1\}}
\newcommand{\syp}{\sigma_P}
\newcommand{\upn}{\frac{1}{(2\pi)^n}}
\newcommand{\equiL}{\underset{L\to \infty}{\sim}}

\newcommand{\tol}{\underset{l\to \infty}{\to}}
\newcommand{\toL}{\underset{L\to \infty}{\to}}

\begin{document}

%%%%% To ease editing, add:

\baselineskip=17pt

%%%%%%%%%%%%%%%%

%% In the running head, give an abbreviation of the title. 
%\titlerunning{Random eigenfunctions}

\title{Betti numbers of random nodal sets  \\
of elliptic pseudo-differential operators}

\author{Damien Gayet
\and 
Jean-Yves Welschinger}

\date{\today}

\maketitle
%34L20   	Asymptotic distribution of eigenvalues, asymptotic theory of eigenfunctions
%60B20 : random matrices
%58J40
%Pseudodifferential and Fourier integral operators on manifolds
 %	58J50   	Spectral problems; spectral geometry; scattering theory
 %60D05 Geometric probability and stochastic geometry

%%%%%%%%

\begin{abstract}
Given an elliptic self-adjoint pseudo-differential operator $P$ bounded from below, acting 
on the sections of a Riemannian line bundle over 
a smooth closed manifold $M$ equipped with some Lebesgue measure,
we estimate from above, as $L$ grows to infinity, 
the Betti numbers of the vanishing locus of a random section taken in the direct sum of the 
eigenspaces of $P$ with eigenvalues below $L$. 
These upper estimates follow from some equidistribution of the critical points
of the restriction of a fixed Morse function to this vanishing locus. We then consider the examples of the Laplace-Beltrami 
and the Dirichlet-to-Neumann operators associated to 
some Riemannian metric on $M$.\\

%% Keywords are optional
Keywords: {Pseudo-differential operator,  random nodal sets, random matrix.}\\

\textsc{Mathematics subject classification 2010}: Primary 34L20, 58J40 ; Secondary 60D05,  60B20.
\end{abstract}

\section*{Introduction}
Let $M$ be a smooth closed manifold of positive dimension $n$, 
by which we mean a smooth compact $n$-dimensional manifold without boundary.
Let $|dy|$ be a Lebesgue measure on $M$, that is locally the absolute value of some volume form. 
Let $E$ be a real line bundle over $M$ equipped with some
Riemannian metric $h_E$. The space $\Gamma(M,E)$ of smooth  global sections of $E$ inherits from $|dy|$
and $h_E$ the $L^2$-scalar product 
\begin{equation}\label{<>}
 (s,t)\in \gme^2 \mapsto \cg s,t\cd = \int_M h_E(s(y), t(y)) |dy|\in \R.
 \end{equation}
 Let then $P : \gme \to \gme$ be an elliptic pseudo-differential operator of order $m>0$ which is self-adjoint with
 respect to (\ref{<>})
 and bounded from below, see \S \ref{paragraphe II.1}. For every $L\in \R$, 
 we denote by 
 \begin{equation}\label{ul}
 U_L = \bigoplus_{\lambda \leq L} \ker (P - \lambda Id)
 \end{equation}
 and  by $N_L$ its dimension. It is equipped with the restriction 
 $\cg \, , \, \cd_L$ of  (\ref{<>})
 and thus with the associated Gaussian measure $\mu_L$ whose density with respect to the Lebesgue measure
  $|ds|$ of $U_L$ reads at every $s\in U_L$,
 \beq 
 d\mu_L (s)= \frac{1}{\sqrt \pi^{N_L}} e^{-\cg s,s\cd} |ds|.
 \eeq
 What is the expected topology of the vanishing locus $\suz\subset M$ 
 of a section $s$ taken at random in $(U_L, \mu_L ) $? A famous theorem of Courant \cite{CourantHilbert}
 bounds from above the number of connected components of $\suz$ whatever $s$ is, when $P$
 is the Laplace-Beltrami operator associated to a Riemannian metric on $M$.
 In the latter case, the expected value for this number of connected components
 for the round two-sphere has been estimated by F. Nazarov and M. Sodin 
 \cite{NazarovSodin},
 a work partially extended in several directions (see \cite{LerarioLundberg}, \cite{Nicolaescu}, \cite{Sodin}, \cite{SarnakWigman} ). 
 We studied a similar question in real algebraic geometry, where 
 $M$ is replaced by a real projective manifold $X$ and $U_L$ 
 by the space $\R H^0(X, E\otimes L^d)$ of real holomorphic sections of 
 the tensor product of some  holomorphic vector bundle $E$ with 
some ample real line bundle $L$ 
 over $X$ (see  \cite{GaWe2}, \cite{GaWe3}, \cite{GaWe4}, \cite{GaWe5}). We there could estimate from above and below the expected value
 of each  Betti number of $\suz$.
 Our aim now  is, likewise,  to estimate from above the 
 mathematical 
 expectations of all Betti numbers
 of $\suz$ for a random section $s\in U_L$,
 as $L$ grows to infinity, see Corollary \ref{Corollaire 0}.
 This turns out to involve asymptotic estimates of the derivatives of the 
 Schwartz kernel associated to the orthogonal projection onto $U_L$
 which we establish in Appendix \ref{A1}, 
 see Theorem \ref{Theoreme 3}. 
 The asymptotic value of this kernel has been 
 computed by L. H\"ormander in \cite{Hormander},
 after Carleman \cite{Carleman} and G\"arding \cite{Garding}
 and for some derivatives, it is given by Safarov and Vassiliev in 
 \cite{SaVa}, but we could not find a general result for all derivatives in the literature.

 Let us now formulate our main result.
When $n\geq 2$, we choose a Morse function $p : M \to \R$ and set 
 $$\Delta_L = \{s \in U_L \, | \, s \text{ does not vanish transversally
  or } p_{|\suz} \text{ is not Morse} \}.$$
Then, for every $s\in U_L\setminus \Delta_L$ and every $i\in \{0, \cdots, n-1\}$, 
we introduce the empirical measure 
$$ \nu_i( s)= \sum_{x\in Crit_i (p_{|\suz})\setminus Crit(p)} \delta_x,$$
where $Crit(p)$ denotes the critical locus of $p$, $Crit_i (p_{|\suz})$
the set of critical points of index $i$ of $p_{|\suz}$
and $\delta_x$ the Dirac measure at $x$. 
When $n=1$, we set
$$ \nu_0( s)= \sum_{x\in \suz} \delta_x.$$
The mathematical expectation
of $\nu_i$ is defined as 
$$\mathbb E(\nu_i)= \int_{U_L \setminus \Delta_L}
\nu_i (s) d\mu_L (s).$$

Recall that the pseudo-differential operator $P$
has a (homogenized) principal symbol $\si_P : T^*M \to \R$
which is homogeneous of degree $m$, see 
Definition \ref{def 2.2},
and we set 
\begin{equation}\label{KP}
K = \{\xi \in T^*M \, | \,  \si_P (\xi)\leq 1\}.
\end{equation}
The volume of $K$ for the Lebesgue measure $|d\xi|$
induced on the fibres of $T^*M$ by $|dy|$ is encoded by
the function 
\begin{equation}\label{c0}
c_0 : x\in M \mapsto \frac{1}{(2\pi)^n } \int_{K\cap T_x^*M} |d\xi|
\in \R_+.
\end{equation}
It turns out that $K$ together with $|d\xi|$ 
induce a Riemannian metric on $M$, 
namely 
\begin{equation}\label{gp}
 g_P : (u,v)\in T_x M \mapsto \frac{1}{(2\pi)^n} \int_{K\cap T_x^*M} \xi(u)\xi(v) |d\xi|
 \end{equation}
and we denote by $|dvol_{P}|$ the associated Lebesgue
measure of $M$.

\begin{thm}\label{Theoreme 0.1}
Let $M$ be a smooth closed manifold of dimension $n$ equipped with a Morse
function $p$ and a Lebesgue measure $|dy|$. Let $(E,h_E)$
be a Riemannian real line bundle over $M$ 
and 
$P : \gme \to \gme$ be an 
elliptic self-ajdoint pseudo-differential operator of order $m>0$ 
which is bounded from below. Then,  for every $\izn$, 
\begin{equation}\label{c0equation}
 \frac{1}{L^{\frac{n}{m}}}\mathbb E(\nu_i)\underset{L\to \infty}{\to} \frac{1}{\sqrt \pi^{n+1}\sqrt c_0}
\mathbb E(i,\ker dp) |dvol_{P}|.
\end{equation}
\end{thm}
The convergence given by (\ref{c0equation})  is the weak convergence
on the whole $M$. 
Also, in Theorem \ref{Theoreme 0.1}, $\mathbb E(i,\ker dp)$ denotes,
for every point $x\in M$,
 the expected determinant of random symmetric operators
of signature $(i,n-1-i)$ on $\ker d_{|x}p$
when $n>1$, see (\ref{ek}),
while  it equals $1$ when $n=1$. Namely, 
 $P$ together with $|dy|$ induce a Riemannian metric 
 $\cg \ , \cd_P$
 on the space
$Sym^2(TM)$ of symmetric bilinear forms on $T^*M$, 
which reads for every $(b_1, b_2)\in Sym^2(TM)^2$,
\begin{equation}\label{sym}
\cg b_1, b_2\cd_P = 
\frac{1}{(2\pi)^n} \Big(
\int_{K} b_1(\xi)b_2 (\xi)|d\xi| -
\frac{1}{\int_{K} |d\xi|} 
\iint_{K^2} b_1(\xi) b_2(\xi')|d\xi||d\xi'|\Big),
\end{equation}
where in the right-hand side of (\ref{sym})
the quadratic forms associated to $b_1$ and $b_2$
are also denoted by $b_1$ and $b_2$, by
abuse of notation. 
The first term in the right-hand side of (\ref{sym})
already defines a natural Riemannian metric
on $Sym^2 (TM)$, see \S \ref{tensoral},
but the one playing a r\^ole in Theorem \ref{Theoreme 0.1}
is indeed (\ref{sym}), where the second 
term induces some correlations similar
to the ones already observed by L. Nicolaescu in \cite{Nicolaescu}.
By duality and restriction to $(\ker dp)^*$, (\ref{sym}) induces a Riemannian metric on 
$Sym^2((\ker dp)^*)$ 
see \S \ref{para 3.3.1},
with Gaussian measure $\mu_P$.  Let $Sym^2_i((\ker dp)^*)$
be the open cone of non-degenerated symmetric bilinear forms
of index $i$ on $\ker dp$. 
We set 
\begin{equation}\label{ek}
 \mathbb E(i,\ker dp) = \int_{Sym^2_i((\ker dp)^*)}
|\det \beta| d\mu_P(\beta),
\end{equation}
where $\det \beta$ is 
 computed with respect to the metric
$g_P$ restricted to $\ker dp$ and given by (\ref{gp}).

From Theorem \ref{Theoreme 0.1} we thus know
that the critical points of index $i$ of $p_{|\suz}$ 
 equidistribute in the manifold $M$ with respect to 
  $g_P$, with a density involving random symmetric
endomorphisms of $\ker dp\subset TM$.
Let us mention two consequences of Theorem \ref{Theoreme 0.1}.
First, for every $s\in U_L\setminus \Delta_L$,
we denote by $m_i(s) $ the $i$-th Morse number of $\suz$,
that is $$m_i(s)= \inf_{f \text{ Morse on }\suz} \# Crit_i(f) $$
and set 
\begin{equation}\label{emi}
\mathbb E(m_i)= \int_{U_L\setminus \Delta_L}
m_i(s) d\mu_L (s).
\end{equation}
From Morse theory we know that these Morse numbers bound from above
all $i$-th  Betti numbers $b_i$ of $\suz$, whatever the coefficient rings are.
% % % % % % % % % % % % % % %
\begin{cor}\label{Corollaire 0}
Under the hypotheses of Theorem \ref{Theoreme 0.1}, when $n\geq 2$,
$$\limsup_{L\to \infty} \frac{1}{L^{\frac{n}{m}}} \mathbb E(m_i)
\leq 
\frac{1}{\sqrt \pi^{n+1}}
\inf_{p \text{ Morse function on } M} \int_M \frac{1}{\sqrt c_0}\mathbb E(i,\ker dp)|dvol_{P}|,
$$
while when $n=1$, we have the convergence
$$ \frac{1}{L^{\frac{1}{m}}} \mathbb E(b_0) \underset{L\to \infty}{\to}
\frac{1}{\pi}
\int_M \frac{1}{\sqrt c_0}|dvol_{P}|.
$$

\end{cor}
Theorem \ref{Theoreme 0.1}
also specializes to the case of the Laplace-Beltrami
operator $\Delta_g$ associated to some  Riemannian metric $g$ 
on $M$. In this case, we denote by $|dvol_g|$ the Lebesgue 
measure associated to $g$ and by $Vol_g(M)$
its total volume $\int_M |dvol_g|$.
\begin{cor}\label{Corollaire 0'}
Let $(M,g)$ be a closed Riemannian manifold of positive dimension $n$
equipped with a Morse function $p : M \to \R$. 
Then, when $n\geq 2$, for every $\izn$, 
$$\frac{1}{\sqrt L^n}
\mathbb E(\nu_i)\toL \frac{\mathbb E(i,n-1-i)}{\sqrt \pi^{n+1}\sqrt {(n+2)(n+4)^{n-1}}}
 |dvol_g|,$$
where the convergence is weak on $M$.
In particular, 
$$\limsup_{L\to \infty} \frac{1}{\sqrt L^n} \mathbb E(m_i) \leq 
\frac{\mathbb E(i,n-1-i)}{\sqrt \pi^{n+1}\sqrt{ (n+2)(n+4)^{n-1}}}
 Vol_g(M).
$$
When $n=1$, 
$\frac{1}{\sqrt L}
\mathbb E(\nu_0)\toL \frac{1}{\pi\sqrt 3} |dvol_g|$
so that $\frac{1}{\sqrt L} \mathbb E(b_0)\toL \frac{1}{\pi\sqrt 3} Vol_g(M).$
\end{cor}
The  case  $n=1$ in Corollary \ref{Corollaire 0'} turns out also to follow from the volume computations
carried out by P. B\'erard in \cite{Berard}.
Note that in Corollary \ref{Corollaire 0'}, $\mathbb E(\nu_i)$ 
is defined using 
$P= \Delta_g$ as a differential operator, 
so that $m=2$ with the notations of Theorem \ref{Theoreme 0.1}. 
Moreover,
\begin{equation}\label{ein}
 \mathbb E(i,n-1-i)= \int_{Sym(i,n-1-i, \R)} |\det A | d\mu(A),
\end{equation}
where $Sym(i,n-1-i, \R)$ denotes  the open cone of non degenerated
symmetric matrices of index $i$, size $(n-1)\times (n-1)$ 
and real coefficients, while $\mu$
denotes the Gaussian measure on $Sym(n-1,\R)$ associated to
the scalar product 
\begin{equation}\label{scal3}
 (A,B)\in Sym(n-1,\R)^2 \mapsto \frac{1}{2}\Tr (AB)+ \frac{1}{6}(\Tr A)(\Tr B)\in \R,
 \end{equation}
see \S \ref{Laplace}.
This measure differs from the standard
GOE measure on $Sym(n-1, \R)$. 
When $M$ is a surface for example, Corollary \ref{Corollaire 0'}
implies that for $i\in \{0,1\}$,
$$ \limsup_{L\to \infty}\frac{1}{L} \mathbb E(m_i)\leq \frac{1}{8\pi^2} Vol_g(M).$$

For large values of the dimension $n$,
we observe some exponential decrease of the upper estimates
given by Corollary \ref{Corollaire 0'} away from the mid-dimensional Betti numbers. This exponential
decrease given by 
Proposition \ref{prop} is similar to the one given by 
Theorem 1.6 of \cite{GaWe4}. 

\begin{Proposition} \label{prop}
For every $\epsilon>0$,  there exist $\delta>0$ and $C>0$
such that for every  smooth closed Riemannian manifold $M$ 
of positive dimension $n$,
$$\limsup_{L\to \infty} \frac{1}{N_L}\sum_{|\frac{i}{n}-\frac{1}{2}|\geq \epsilon} \mathbb E(m_i)\leq C\exp (-\delta n^2).$$
\end{Proposition}
In particular, $$\limsup_{L\to \infty}\frac{1}{N_L}\mathbb E (b_0) \to_{n\to \infty}
0,$$
which has to be compared with the Courant 
upper bound $b_0\leq N_L$, see \cite{CourantHilbert}.
Again, in Proposition \ref{prop}, $\mathbb E(m_i)
$ is defined using $P= \Delta_g$ a s a differential operator. 

As a second example, Theorem \ref{Theoreme 0.1}
specializes to the case of the Dirichlet-to-Neumann
operator on the boundary $M$ of some compact
Riemannian manifold $(W,g)$, see \S \ref{Dirichlet}. 
We then obtain 
\begin{cor}\label{d2n}
Let $(W,g)$ be a smooth compact Riemannian manifold 
of positive dimension $n+1$ with boundary $M$,
$\Lambda_g$ be the Dirichlet-to-Neumann
operator on $ M$, and 
$p :  M \to \R$ be a fixed Morse function.
Then, when $n\geq 2$, for every $\izn$, 
$$\frac{1}{L^n}
\mathbb E(\nu_i)\toL \frac{\mathbb E(i,n-1-i)}{\sqrt \pi^{n+1}\sqrt {(n+2)(n+4)^{n-1}}}
 |dvol_g|,$$
where the convergence is weak on $ M$
and $|dvol_g|$ is the volume form on $M$
induced by $g$.
In particular, 
$$\limsup_{L\to \infty} \frac{1}{L^n} \mathbb E(m_i) \leq 
\frac{\mathbb E(i,n-1-i)}{\sqrt \pi^{n+1}\sqrt{ (n+2)(n+4)^{n-1}}}
 Vol_g(M).
$$
When $n=1$, 
$\frac{1}{L}
\mathbb E(\nu)\toL\frac{1}{\pi\sqrt 3} |dvol_g|$
so that $\frac{1}{L} \mathbb E(b_0)\toL \frac{1}{\pi\sqrt 3} Vol_g( M).$
\end{cor}

In the  first section
we study the general case of an ample finite dimensional 
subspace $U$ of $\gme$ equipped with any scalar product, 
see Definition \ref{definition 1}.
We get estimates similar to the ones
given by Theorem \ref{Theoreme 0.1} and Corollary \ref{Corollaire 0},
in terms of the Schwartz kernel associated to $U$ and its derivatives, see 
\S \ref{paragraphe 1.4}.
The second section is devoted to 
the special case of $U=U_L$ for 
some elliptic self-adjoint pseudo-differential
operator bounded from below, see (\ref{ul}), and to the proofs of 
Theorem \ref{Theoreme 0.1} and Corollary \ref{Corollaire 0}.
The third section is devoted to examples, namely the case
of the Laplace-Beltrami and Dirichlet-to-Neumann operators, where we prove
 Corollary \ref{Corollaire 0'}, Proposition \ref{prop}, and Corollary \ref{d2n}.
In the last section we 
discuss some related problems 
which we plan to consider in a separated paper. 
We finally give in Appendix  \ref{AppB} 
 several auxiliary results, in particular 
 the proof of Theorem \ref{Theoreme 3}, which provides
 estimates of the derivatives
 of the Schwartz kernel associated to $U_L$. 

\textit{Aknowledgements.} The research leading to these results has received funding
from the European Community's Seventh Framework Progamme 
([FP7/2007-2013] [FP7/2007-2011]) under
grant agreement $\text{n}\textsuperscript{o}$ [258204].
We are grateful to Yves Colin de Verdi\`ere for suggesting 
the example given in \S \ref{Dirichlet}  and to 
Antonio Lerario for pointing out 
the reference \cite{Nicolaescu} to us.

\tableofcontents

%%%%%%%%%%%%%%%%%%%%%%%%%%%%%%%%%%%%%%%%
\section{Morse numbers 
of the vanishing locus of 
random sections}\label{paragraphe 1}
Let $M$ be a smooth manifold of positive dimension $n$,
$E\to M $ be a real line bundle 
and $p: M\to \R$ be a Morse function.
We denote by $\bh$ the singular
foliation by  level sets of $p$
and for every $ x\in M\setminus {Crit(p)}$ we set
$$H_x = T_x \bh = \ker \dxp.$$
\subsection{Ample linear subspaces and incidence varieties}\label{paragraphe 1.1}

For every $l\geq 0$, we denote by $\jle$
the fibre bundle of $l$-jets of sections of $E$
and for every $m\geq l\geq 0$, we denote
by $\pi^{m,l} : \jme \to \jle$
the tautological projections which restricts the $m$-jets to $l$-jets. The jet maps
are denoted by
$$ j^l : s\in \gme \mapsto j^l(s)\in \Gamma(M, \jle).$$
Recall that the kernel of $\pi^{l+1,l} $
is canonically isomorphic to the bundle 
$Sym^{l+1}(T^*M)\otimes E$ of 
symmetric $(l+1)$-linear forms 
on $TM$ with values in $E$. 
In particular, any Riemannian metric on $\jle$ induces an isomorphism
$$\jle \cong S^l (T^*M)\otimes E,$$
where $S^l (T^* M) = \bigoplus_{k=0}^l Sym^k(T^*M).$

Let $U\subset \gme$ be a linear subspace of positive dimension $N$ 
and $ \us = M\times U $ be 
the associated  rank $N$ trivial bundle over $M$. 
The maps $j^l $ define bundle morphisms
$$ j^l : (x,s) \in \us \mapsto (x,j^l (s)_{|x})\in \bj^l(E).$$
\begin{defin}\label{definition 1}(compare Def. 2.1 of \cite{Nicolaescu})
The vector subspace $U$ of $\gme$ is said to be $l$-ample
if and only if  the morphism 
$j^l : \us \to \bj^l(E)$
is onto. It is said to be ample if and only if  it is $1$-ample.
\end{defin}

We also need a relative version of this ampleness property. 
For every $l\geq 0$, we denote by $\jleh\to \mcp$ 
the fibre bundle of $l$-jets of restrictions
of sections of $E$ to the leaves of $\bh$.
If $x\in \mcp$ and $\bh_x = p^{-1}(p(x))$,
then the fibre of $\jleh$ over $x$ 
is the space of $l$-jets at $x$ of sections of the restriction $E_{|\bh_x}.$
These bundles are likewise equipped with projections 
$$ \pi^{m,l} : \jmeh \to \jleh,$$
$m\geq l\geq 0$
and with jet maps 
$$ j^l_\bh : s\in \gme \mapsto  j^l_\bh (s)\in 
\Gamma(\mcp, \jleh).$$ 
These jet maps induce bundle morphisms 
\begin{equation}\label{jlh}
j^l_\bh : (x,s) \in \us_{|M\setminus Crit(p)} \mapsto
(x,j^l_\bh (s)_{|x})\in \jleh.
\end{equation}
\begin{defin}\label{definition 2.2}
The linear subspace $U$ of $\gme$ is said to be relatively $l$-ample
if and only if  the bundle morphism 
$j_\bh^l :  \us_{|M\setminus Crit(p)} \to
 \jleh$
is onto. It is said to be relatively ample if and only if  it is relatively  $1$-ample.
The kernel of $j^l_\bh$ is then called the $l$-th incidence variety
and denoted by $\bi^l$. 
\end{defin}
The incidence varieties given by Definition \ref{definition 2.2} are equipped with projections 
\beq
\pi_M : &(x,s)\in \bi^l &\mapsto x\in  \mcp \text{ and }\\
 \pi_U : &(x,s)\in \bi^l &\mapsto s \in U,
 \eeq
see \ref{A0} for further properties. 
We set
\begin{eqnarray}\label{Delta1} 
\Delta_0 &=& \{s\in U\, | \ s \text{ does not vanish transversally} \} \text{ and if } n\geq 2, \nonumber\\
\Delta_1 &=& \Delta_0 \cup \{s\in U\setminus \Delta_0\ | \ p_{|s^{-1}(0)} \text{ is not Morse}. \}
\end{eqnarray}
Then, for every $i\in \{0, \cdots, n-1\}$, we set
$$ \bi_i^1 = \{\xs \in (\mcp) \times (U\setminus \Delta_1) \ | \ s(x)= 0 
\text{ and } x\in Crit_i (p_{|\suz}) \},$$
where $Crit_i (p_{|\suz}) $ denotes the set of critical points
of index $i$ of the restriction of $p$ to $\suz$. 
The disjoint union 
$ \bi^1_0\cup \cdots \cup \bi^1_{n-1}$ provides
a partition of $\bi^1\setminus \pi^{-1}_U (\Delta_1),$
see Appendix \ref{A0}.

These incidence varieties
equip $\us_{|\mcp}$ with some filtration 
whose first graded maps read 
\beq
gr^0 : (x,s_0)\in  \us/\bi^0 &\mapsto & s_0(x) \in E 
\eeq
and 
\beq
gr^1 : (x,s_0, s_1)\in  \us/\bi^0 \oplus  \bi^0/\bi^1 &\mapsto &(s_0(x), \nabla s_{1|H_x})\in  E \oplus (H^*\otimes E ).
\eeq
Finally, we set 
\beq
H^\circ &=& \{\lambda \in T^*M \, | \, \lambda _{|H}= 0\} \text{ and }\\
 j: (x,s)\in \bi^1 &\mapsto&  (x, \nabla s , \nabla^2 s _{|H_x})\in (H^\circ \oplus Sym^2(H^*))\otimes E
\eeq
when $n\geq 2$,
while we set 
\beq j_0: (x,s) \in \bi^0 &\mapsto & (x, \nabla s) \in  T^*M \otimes E
\eeq
when $n=1$.
Note that $\det (gr^1) = \det (j^1_\bh) : \det (\us /\bi^1)\to \det (H^*)\otimes (\det E)^n$
and that for every $(x,s)\in \bi^1$, $j(x,s)$ induces the morphisms 
\beq 
j(x,s) : 
T_xM /H_x \oplus H_x &\to& E_x \oplus (H^*_x \otimes E_x) \text { and }\\
\det (j(x,s)) : \det (T_xM)&\to& \det (H^*)\otimes (\det E)^n.
\eeq

\subsection{The induced Riemannian metrics}\label{paragraphe 3}
% % % % % % % % % % % % % % % % % %
\begin{lem}\label{Lemme 3}
Let $F$, $G$ be two finite dimensional real vector spaces  and
$A : F\to G$ be an onto linear map. Let $\langle, \rangle_F$
be a scalar product on $F$ and $\# : F^* \to
F $ be the associated isomorphism. Then, the composition 
$ (A \# A^*)^{-1}: G\to G^*$
defines a scalar product $\langle \, , \, \rangle_G$ on $G$.
Moreover, if $\mu_F$ (resp. $\mu_G$)
denotes the Gaussian measure associated to $\langle \, , \, \rangle_F$
(resp. $\langle\, , \, \rangle_G$), then $\mu_G = A_*\mu_F$.
\end{lem}
Let $|df|$ (resp. $|dg|$) be the Lebesgue measure 
associated to $\langle \, , \, \rangle_F$ (resp. $\langle\, , \, \rangle_G$).
Then $$d\mu_F(f) = \frac{1}{\sqrt \pi^{\dim F}}e^{-\| f\|^2} |df|$$
and $d\mu_G(g) = \frac{1}{\sqrt \pi^{\dim G}}e^{-\| g\|^2} |dg|,$
where $\| f\|^2 = \langle f, f\rangle_F$ and
$\| g\|^2 = \langle g, g\rangle_G$. 
\bpr
Let $g^*_1$, $g_2^*\in G^*$. Then
$
\cg g_1^*, g_2^* \cd_{G^*} = g_2^* (A \# A^* (g_1^*))
= A^* (g_2^*)(\#  A^*(g_1^*))
  = \cg \# A^* (g_2^*), \# A^*(g_1^*)\cd_F.$
Since $A^*$ is injective, we deduce that $\cg \, , \, \cd_{G^*}$
is a scalar product on $G^*$ and hence that $\cg \, , \, \cd_{G}$
is a scalar product on $G$.
Moreover, $\#  A^* : G^* \to (\ker A)^\perp $ 
is an isometry, so that 
$A : (\ker A)^\perp \to G$ is an isometry. Since $\mu_F$ 
is a product measure, we deduce that $\mu_G = A_* \mu_F$.
\epr
\begin{defin}
Under the hypotheses of Lemma \ref{Lemme 3}, 
$\cg \, , \, \cd_{G}$ (resp. $\mu_G$ )
is called the push-forward of $\cg \, , \, \cd_{F}$
(resp. $\mu_F$) under $A$.
\end{defin}

\begin{defin}\label{defpb}
Let $U\subset \gme$ be an ample finite dimensional linear 
subspace, which
is equipped with a scalar product $ \cg\ , \cd$. 
The latter induces a Riemannian metric on the trivial
bundle $\us$ which restricts to a metric on $\mathcal I^l$, $l\in \Nn$.
We denote by $\mu_{\bi^l}$ the associated
Gaussian measure and by 
\begin{itemize}
\item $g^1 $ the push-forward on $E\oplus (H^*\otimes E)$  of  $ \cg\ , \cd$ under  $gr^1$,
\item $h^l$ the push-forward on $\bj^l(E_{|\bh}) $  of  $ \cg\ , \cd$ under  $j^l_\bh$  
and 
\item  $h$ the push-forward  on 
$Im(j)  \subset  (H^\circ \oplus Sym^2(H^*))\otimes E$ of  $ \cg\ , \cd$ under $j$,
\end{itemize}
see \S \ref{paragraphe 1.1} and Lemma \ref{Lemme 3}.

When $n=1$, we denote by 
\begin{itemize}
\item $g^0 $ the push-forward on $E$  of  $ \cg\ , \cd$ under  $gr^0$,
\item  $h_0$ the push-forward  on 
$Im(j_0)  \subset  T^*M\otimes E$ of  $ \cg\ , \cd$ under $j_0$.
\end{itemize}
 \end{defin}
\begin{defin}\label{Definition 2} The Schwartz kernel of $(U, \cg \, , \, \cd)$
is the section
$e$ of $\us \otimes E$ satisfying for every 
$ s\in U$ and $x\in M, \ s(x)= \cg e_x, s\cd.$
\end{defin}
Note that if $(s_1, \cdots, s_N)$ denotes an orthonormal basis of $U$,
then for every $x\in M$, $ e_x = \sum_{i=1}^N s_i(x)s_i.$
The metrics $g^1$, $h^l$ and $h$ given by Definition \ref{defpb} can be computed 
in terms of the Schwartz kernel $e$, 
as follows from Lemma \ref{Lemme 4} and \ref{Lemme 5},
compare \cite{DSZ1}, \cite{Nicolaescu}
% % % % % % % % % % % % % %
\begin{lem}\label{Lemme 4}
Let $E$ be  a real line bundle  over a smooth manifold $M$
equipped with a Morse function. Let $U$ be a 
finite dimensional  linear subspace 
of $\gme$ which is relatively $l$-ample for  $l\in \Nn^*$ and equipped with a scalar product. 
Let $e$ be its Schwartz kernel. Then, the metrics $h^l$ and $g^1$
are given by the restriction to the diagonal of 
$ (j^l_\bh j^l_\bh e)^{-1} $ and $(gr^1 gr^1 e)^{-1}.$
\end{lem}
% % % % % % % % % % % % % % %
Note that $e$ is a section of $E\boxtimes E$ over $M\times M$, 
so that $j_\bh^l j_\bh^l e $ (resp. $gr^1 gr^1 e$),
which applies $j^l_\bh$ (resp. $gr^1$) on each variable 
of $e$, is a section of
$ \bj^l(E_{|\bh})^{\boxtimes 2}$ (resp. $(E\oplus (H^*\otimes E))^{\boxtimes 2}$). Its restriction to the diagonal thus
defines
a symmetric  bilinear form on $\bj^l(E_{| \bh})^*$ (resp. $(E\oplus (H^*\otimes E))^*$).
\bpr Let $\theta^*\in \bjlh^*$ and $s\in U$. Then, $s=\cg e,s\cd$
and $(j^{l*}_\bh \theta^* )(s)= \cg \theta^*(j_\bh^l e), s\cd$. 
Consequently, $ \# (j^l_\bh)^* \theta^* = \theta^* (j^l_\bh e)$
and 
$ j^l_\bh \# (j^l_\bh)^* = j_\bh^l j^l_\bh e$. 
Likewise, $gr^1 \# gr^{1*} = gr^1 gr^1 e.$
\epr
\begin{lem}\label{Lemme 5}(Compare appendix A of \cite{Nicolaescu})
Let 
$A : F\to G$
be a linear map between two 
real finite dimensional vector spaces. Let $K_F$ (resp. $K_G$)
be a subspace of $F$ (resp. $G$) such that $A(K_F)\subset K_G$
and $a : K_F \to K_G$ be the restriction of $A$. 
Let $\cg \, , \, \cd_F$ 
be a scalar product on $F$ and let $K_F$ be equipped 
 with its restriction. 
Let $L_G$ be a complement subspace of $K_G$ in $G$ and 
$ b : K_F^\perp \to K_G$ (resp. $c: K_F^\perp \to L_G$)
be such that 
$$ A =\begin{bmatrix}
a & b\\ 0 & c
\end{bmatrix} : K_F\oplus K_F^\perp \to K_G\oplus L_G.$$
Then,  
$$A\# A^*= \begin{bmatrix}
a\#a^* + b\# b^* & b\# c^*\\
 c\#b^* & c\# c^*
 \end{bmatrix}. $$
$ \square$
 \end{lem}
\begin{rem}\label{Remarque a}
Since
 $a\# a^* = (a\#a^* + b\# b^*) - b\# c^*(c\# c^*)^{-1} c\# b^*$, we deduce
 from Lemma \ref{Lemme 5}
 that the scalar product $(a\#a^*)^{-1}$ 
 can be computed  from $(A\#A^*)^{-1}$. 
Applying Lemma \ref{Lemme 5} to 
$$\left\{\begin{array}{lcl}
F&=& \us,\\
G&=& \bj^1 (E)\times_M \bj^2 (E,\bh),\\
K_F &=& \bi^1 \text { and}\\
K_G &=& Im (j)
\subset (H^\circ\oplus Sym^2 (H^*))\otimes E,
\end{array}\right.
$$
we deduce that the metric $h$ can be computed
in terms of the Schwartz kernel $e$ of $U$ 
and the jet maps $j^1$ and $j^2_\bh$. 
\end{rem}

\subsection{Distribution of critical points}\label{paragraphe 1.4}
\subsubsection{The main result}
Let $U\subset \gme $ be a relatively $l$-ample linear subspace of finite dimension $N$, see
Definition \ref{definition 2.2}.
We equip $U$ with a scalar product $\cg\ , \ \cd$
and denote by $\mu_U$ the associated Gaussian measure, 
so that at every point $s\in U$  its density against the Lebesgue measure $|ds|$ on $U$
equals $\frac{1}{\sqrt \pi^N}e^{-\cg s,s\cd}.$
Then, for every $i\in \{0, \cdots n-1\}$ and every $s\in U\setminus \Delta_1$, where $\Delta_1$ is given by (\ref{Delta1})
we set 
\beq
 \nu_i (s)&=& \sum_{x\in Crit_i (p_{|\suz})\setminus Crit(p)} \delta_x \\
\mathbb E(\nu_i)&=& \int_{U\setminus \Delta_1} \nu_i(s)d\mu_U(s)
\eeq
when $n\geq 2$, while when $n=1$, we set 
\beq  
\nu_0 (s)&=& \sum_{x\in \suz} \delta_x\\
\mathbb E(\nu_0)&=& \int_{s\in U\setminus \Delta_0} \nu_0 (s) d\mu_U(s).
\eeq
Note that  we have no control a priori on the number
of critical points of the restriction of $p$ to $\suz$, 
so that $\mathbb E(\nu_i)$ may not be well defined.
% % % % % % % % % % % % % %
\begin{thm}\label{Theoreme 1}
Let $E$ be a real line bundle over a smooth
$n$-dimensional 
manifold $M$ equipped with a Morse function. Let $U\subset \gme$
be a finite dimensional relatively ample linear subspace equipped with a scalar product. Then, when $n\geq 2$, 
for every $i\in \{0, \cdots, n-1\},$
\begin{equation}\label{eni}
\mathbb E(\nu_i)= \frac{1}{\sqrt \pi^n} 
\iint_{(H^\circ \oplus Sym^2_i(H^*))\otimes E}
|(\alpha, \beta)^* dvol_{g^1}| j_*d\mu_{\bi^1}(\alpha,\beta).
\end{equation}
Moreover, this measure has no atom and its density 
with respect to any 
 Lebesgue measure lies in $C^\infty (M\setminus Crit(p))$. 
If in addition at every point $x\in Crit(p)$ 
the jet map $j^1 : U\to \bj^1(E)_{|x}$ is onto,
then this density lies in 
$L^1_{loc} (M)$, so that $\mathbb E(\nu_i)$
defines a measure on the whole $M$.
When $n=1$, 
$$ \mathbb E(\nu_0)= \frac{1}{\sqrt \pi} 
\int_{T^*M\otimes E}
|\alpha^* dvol_{g^0}| j_{0_*}d\mu_{\bi^0}(\alpha).
$$
\end{thm}
% % % % % % % % % % % %
Theorem \ref{Theoreme 1} describes the expected distribution 
of critical points of the restriction $p_{|\suz}$. 
Every pair $(\alpha,\beta)\in  (H^\circ\oplus Sym^2 (H^*))\otimes E$
defines a morphism
$$ (\alpha, \beta): (TM/H )\oplus H \to E\oplus (H^*\otimes E),$$
while the bundle $E\oplus (H^*\otimes E)$ is equipped
with the metric $g^1$ 
and its associated volume form $dvol_{g^1}$, see Definition \ref{defpb}.
It follows that  $(TM/H) \oplus H$ inherits the $n$-form 
$(\alpha, \beta)^*dvol_{g^1}$. The latter
induces a $n$-form on $TM$, also denoted by 
$(\alpha, \beta)^*dvol_{g^1}$, 
since $\det (TM)$ is canonically isomorphic to $\det ((TM/H)\oplus H).$
Finally, we have denoted by $Sym^2_i(H^*)$ the open cone of 
non-degenerated symmetric bilinear  forms  of index $i$ on $H$.
 Recall that the index of a symmetric bilinear
 form is the maximal dimension of a subspace 
on which the form restricts to a negative definite one.
Note that the form $(\alpha, \beta)^*dvol_{g^1}$ depends
polynomially on $(\alpha, \beta)$, so that it
is integrable with respect to the Gaussian measure $j_* \mu_{\bi^1}.$
Note finally that from Lemma \ref{Lemme 4}
and Remark \ref{Remarque a}, both $g^1$ 
and $j_*d\mu_{\bi^1}$
can be computed in terms of the Schwartz kernel
of $(U, \cg \, , \, \cd)$, see Definition \ref{Definition 2}.
\bpr
By definition, $\mathbb E(\nu_i)= (\pi_{M|\bi^1_i})_*\pi_U^*d\mu_U$
since the measure of $\Delta_1$ vanishes by Lemma \ref{Lemme 1}. 
From the coarea formula, see Theorem 3.2.3 of \cite{Federer} or Theorem 1 of \cite{SS}, we get
$$ (\pi_{M|\bi^1_i})_*\pi_U^* d\mu_U = 
\frac{1}{\sqrt \pi^n}
\int_{\bi_i^1}
| dvol_{((d\pi_M\circ d\pi_U^{-1}) \# (d\pi_M \circ d\pi_U^{-1})^*)^{-1}}|
d\mu_{\bi^1},
$$
Note indeed that $\bi^1$ has codimension $n$
in $\us$, so that the normalization in $d\mu_{\bi^1_i}$
and $d\mu_U$ differs by a factor $1/\sqrt \pi^n$. 
For every $\xs\in \bi^1$, 
$$ T_{\xs} \bi^1 = \{\xsd \in T_{\xs} \us \ | \ j^1_\bh (\sd) + 
\nabla^\bj_{\xd} (j^1_\bh (s))= 0\},
$$
see (\ref{TI1}), so that
$d_{|\xs} \pi_M \circ d_{|\xs} \pi_U^{-1} = -(\nabla^\bj (j^1_\bh (s)))^{-1}\circ j^1_\bh. $
The operator $\nabla^\bj (j^1_\bh(s))$ is invertible since $s\in U_L\setminus \Delta_1$, see Remark \ref{rema}.
It follows that the determinant of the morphism $\us/\bi^1 \to TM$
induced by $d_{|\xs} \pi_M \circ d_{|\xs} \pi_U^{-1}$ coincides with the one
of 
$$-j(s)^{-1}\circ gr^1 : \us/\bi^0 \oplus \bi^0/\bi^1 \to TM/H \oplus H$$
via the canonical isomorphisms 
$\det (\us/\bi^1)\cong \det (\us/\bi^0 \oplus \bi^0/\bi^1)$
and $\det (TM)\cong \det (TM/H \oplus H).$
We deduce that 
$$ dvol_{((d\pi_M\circ d\pi_U^{-1}) \# (d\pi_M \circ d\pi_U^{-1})^*)^{-1}}
= dvol_{((j(s)^{-1}\circ gr^1) \# (j(s)^{-1}\circ gr^1)^*)^{-1}} =
 j(s)^* dvol_{g^1}.$$
Using the substitution $(\alpha, \beta)=j(s)$, we conclude that 
$$ \mathbb E(\nu_i)= \frac{1}{\sqrt \pi^n}
\int_{(H^\circ\oplus Sym^2_i(H^*))\otimes E}
|(\alpha, \beta)^* dvol_{g^1}| j_*\mu_{\bi^1}(\alpha,\beta).
$$
Note that $g^1$ is a smooth metric on
$E\oplus (H^*\otimes E)$ since $\mu_{\bi^1}$
is a smooth family of Gaussian measures on $\bi^1$
and $j$ a smooth morphism. We deduce 
that $\mathbb E(\nu_i)$ has no atom and that its density 
with respect to any Lebesgue measure on $M$ belongs to $C^\infty (M\setminus Crit(p)).$

Now, let us assume in addition that  at every critical point $x$ of $p$,
the jet map  $j^1 : U\to \bj^1(E)_{|x}$ is onto and
let us prove
that this density then also belongs to $L^1_{loc}(M)$, so that 
$\mathbb E(\nu_i)$ extends to a measure without atom on the whole $M$.
We denote by $\pi : P(T^*M)\to M$ the projectivization of the cotangent bundle and by
$\tau \subset \pi^*(T^*M)$ the tautological line bundle over $P(T^*M)$. 
From the inclusion $\tau\otimes \pi^*E \to \pi^* (T^*M\otimes E) $ 
we deduce the short exact sequence
$$ 0 \to \tau \otimes \pi^*E \to \pi^* \bj^1(E)\to \pi^*\bj^1(E)/\tau \otimes \pi^* E \to 0.$$
With a slight abuse of notation, we denote by $H\subset \pi^* (TM)$
the codimension one  subbundle 
given by the kernels of the elements of $\tau \setminus \{0\}$
and by $\bj^1(E,H)$ 
the quotient bundle 
$\pi^*\bj^1(E)/\tau \otimes \pi^*(E)$.
Let $V$ be a compact neighbourhood of $Crit(p)$
such that the restriction of the morphism
$ j^1: \us_{|V} \to \bj^1(E)_{|V} $ is onto. 
We deduce a morphism
$j^1 : \pi^* \us \to \pi^*\bj^1(E)$
over $P(T^*M)_{|V}$ which is onto 
and by composition with
the onto map $\pi^* \bj^1(E)\to \bj^1(E,H)$,
an onto morphism $\pi^*{\us} \to \bj^1(E,H)$.
We denote, with an abuse of notation, by $\bi$ the kernel
of the latter and by $g^1$ the metric
that this morphism induces by push-forward on $\bj^1(E,H)$ over $P(T^* M)_{|V}$,
 see Lemma \ref{Lemme 3}.
Now, let $\nabla$ be a torsion-free connection 
on $M$ and let $\nae$ be a connection on $E$.
They define a bundle morphism
$$ \bj : s\in \bi \mapsto (\nabla s_{|H}, \nabla (\nae s)_{|H^2})
\in (\tau \oplus Sym^2 (H^*))\otimes \pi^* E.$$
We then set
$$ \Omega = \frac{1}{\pi^n}\int_\bi 
\bj (s)^* |dvol_{g^1}|
d\mu_\bi (s)=
\frac{1}{\pi^n}
\int_{\tau \otimes \pi^*E}\int_{Sym^2(H^*)\otimes \pi^*E}
(\alpha, \beta)^* |dvol_{g^1}|
(\bj_* d\mu_\bi ) (\alpha, \beta),$$
where $\mu_\bi$ denotes the fiberwise 
Gaussian measure
associated
to the restriction of the metric 
of $\pi^* \us $ to $\bi$. Consequently, 
$\Omega $ provides a section of the fibre bundle $\pi^*\det (T^*M)
$
over the compact $P(T^*M)_{|V}$. 
Let $\omega$ be a volume form on $V$. It
trivializes $\det (T^*M)$ over $V$ and 
$\pi^* \det (T^*M)$ over $P(T^*M)_{|V}$. 
We deduce that there exists a positive constant 
$c>0$ such that $|\Omega| \leq c|\omega|$
over $P(T^*M)_{|V}$. 
However, from Lemma \ref{Lemme 2},
the jet map $j$ on $\bi^1$ 
factors as $j= T\circ \bj$, where 
$T$ denotes the trigonal endomorphism 
of $(H^\circ\oplus Sym^2(H^*))\otimes E$
defined by 
$$ (\alpha, \beta)\mapsto (\alpha, \beta - (\frac{1}{dp}\nabla (dp)_{|H^2})\alpha)$$
and where $\bi^1$ is identified with the pull-back 
$[dp]^*\bi$ by the section $[dp]$ of $P(T^*M)_{|\mcp}
$
defined by the differential of $p$. Finally, 
\beq
\mathbb E(\nu_i) & = & \frac{1}{\pi^n} \iint_{(H^\circ\oplus Sym^2_i(H^*))\otimes \pi^*E}
|T^*(\alpha, \beta)^* dvol_{g^1}|
\bj_* d\mu_{\bi^1} \\
& \leq & 
\frac{C|\omega|}{\pi^n} 
\iint_{(H^\circ\oplus Sym^2_i(H^*))\otimes \pi^*E}
|\det T(\alpha, \beta)|
\bj_* d\mu_{\bi^1} .
\eeq
Since the differential $dp$ vanishes transversally
on $Crit(p)$, 
the function $\det T(\alpha, \beta)$ 
is polynomial in $\alpha, \beta$
and his coefficients are smooth functions on
$\mcp$ with poles of order at most $n-1$
at $\mcp.$ After integration against
the Gaussian measure $\bj_* d\mu_{\bi^1}$, 
we deduce that the function 
$$ 
\int_{H^\circ\otimes E}\int_{Sym^2_i(H^*)\otimes \pi^*E}
|\det T(\alpha, \beta)|
\bj_* d\mu_{\bi^1} $$
is smooth over $\mcp$ with poles of order at most $n-1$
on $Crit(p)$ 
(compare Remark 3.3.3 of \cite{GaWe5}). 
Since $\dim M= n$, 
we deduce that this function belongs to $L^1_{loc} (M)$,
so that $\mathbb E(\nu_i)$ extends to a measure without atom
over the whole $M$.

In the case $n=1$, 
$$ T_{(x,s)}\bi^0 = \{(\xd, \sd)\in T_{(x,s)}\us \, |\, 
\sd (x)+ \nabla_{\xd} s_{|x }= 0\},$$
so that $d_{|\xs} \pi_M \circ d_{|\xs} \pi_U^{-1} = 
-(j_0(s))^{-1}\circ gr^0. $
We deduce that 
$$ dvol_{((d\pi_M\circ d\pi_U^{-1}) \# (d\pi_M \circ d\pi_U^{-1})^*)^{-1}}
= dvol_{((j_0(s)^{-1}\circ gr^0) \# (j_0(s)^{-1}\circ gr^0)^*)^{-1}} =
 j_0(s)^* dvol_{g^0}.$$
Using the substitution $\alpha=j_0(s)$, we conclude that 
$ \mathbb E(\nu)= \frac{1}{\sqrt \pi}
\int_{T^*M\otimes E}
|\alpha^* dvol_{g^0}| j_{0_*}\mu_{\bi^0}(\alpha).
$
\epr
\subsubsection{Mean Morse numbers}\label{paragraphe 1.4.2}
Under the hypotheses of Theorem \ref{Theoreme 1}, 
assume in addition that $M$ is compact without 
boundary. Then, for every $s\in U\setminus \Delta_1$,
$\suz$ is a smooth compact hypersurface of $M$
and for every $i\in \{0, \cdots , n-1\}, $ we set
$$ \mathbb E(m_i)= \int_{U\setminus \Delta_1} m_i (s)d\mu_U(s),$$
see (\ref{emi}).

\begin{cor}\label{Corollaire 1}
Under the hypotheses of Theorem \ref{Theoreme 1}, 
we assume in addition that $M$ is  closed. Then,
for every $\izn$ and every volume form $\omega $ on $M$,
$$ \mathbb E(m_i)\leq \frac{1}{\sqrt \pi^n} \int_M 
\iint_{(H^\circ \oplus Sym^2_i(H^*))\otimes E}
|(\alpha, \beta)^* dvol_{g^1}| j_*d\mu_{\bi^1}(\alpha,\beta).
$$
\end{cor}
\bpr
Corollary \ref{Corollaire 1}
is a consequence of Theorem \ref{Theoreme 1}
after integration of the constant function $1$.
\epr
% % % % % % % % % % %
\subsubsection{An asymptotic result}\label{param}

Let now $(U_L)_{L\in \R_+^*}$ be a family
of finite dimensional linear subspaces of $\gme$
which are ample for $L$ large enough. 
We want to estimate the asymptotic of the
measure $\mathbb E(\nu_i)$ computed by 
Theorem \ref{Theoreme 1} as $L$ grows to infinity.
In order to do so, we need to assume that the 
family $(U_L)_{L\in \R^*_+}$ is tamed in some sense
and from Remark \ref{Remarque a}, we know
that it is sufficient to tame the 
Schwartz kernel $(e_L)_{L\in \R^*_+}$, see
Definition \ref{Definition 2}. However,
we found it convenient to 
tame directly the induced metrics given
by Definition \ref{defpb}, see Definition \ref{tamed}.

\begin{defin}\label{Definition 4}
Let $p,q$ be two positive integers. A
 one-parameter $(p,q)$-group
of endomorphisms of jet bundles is a one-parameter
group $(a_L)_{L\in \R^*_+}$
of diagonalizable endomorphisms on the jet bundles
$\bj^l(E)$, $l\in \Nn$ such that
\begin{enumerate}
\item For every $ 0\leq l\leq m$, the projection
$\pi^{m,l} : \bj^m(E)\to \bj^l(E)$
is $a_L$-equivariant.
\item For every $l\in \Nn$, the restriction of $a_L$
to $\ker \pi^{l+1,l}= Sym^{l+1}(TM^*)\otimes E$ is 
a homothetic transformation of ratio $L^{-p-(l+1)q}.$
\end{enumerate}
\end{defin}
Any such one-parameter $(p,q)$-group of endomorphisms is obtained
in the following way. 
We choose, for every $l\in \Nn$, a complement 
subspace $K_{l+1}$ to $\ker \pi^{l+1,l}$
in $\bj^{l+1}(E)$ and then we require that
$a_L$ preserves $K_{l+1}$ for every $l\in \Nn$, $L\in \R^*_+$. 
The two conditions of Definition \ref{Definition 4}
then determine $(a_L)_{L\in \R^*_+}$
in a unique way.
Note that any metric on $\bj^{l+1}(E)$ provides such 
a complement $K_{l+1}$ to $\ker \pi^{l+1,l}$, 
namely its orthogonal complement
and induces then an isomorphism $\bj^{l+1}(E) \cong S^{l+1}(T^*M\otimes E)$. 
% % % % %
\begin{lem}\label{Lemme 6}
Let $E$ be a real fibre bundle over a smooth  manifold
$M$. Let $(a_L)_{L\in \R^*_+}$ and $(b_L)_{L\in \R^*_+}$
be two one-parameter $(p,q)$-groups of jet bundle
endomorphisms, $p,q>0$. Then, for every $l\in \Nn$,
the composition 
$$ a_L \circ b_L^{-1} : \bj^l(E)\to \bj^l(E)$$ converges
to the identity as $L$ grows to $\infty$.
\end{lem}
% % % % %
\bpr
We proceed by induction on $l\in \Nn$.
When $l=0$, $a_L$ and $b_L$ 
are homothetic transformations of ratio $L^{-p}$ on
$\bj^0(E)$, so that $a_L\circ b_L^{-1}$ 
equals the identity for every $L\in \R^*_+$.
Let us now assume that Lemma \ref{Lemme 6} holds true up to $l\in \Nn$
and prove it for $l+1$.
The endomorphisms $a_L$ and $b_L$ are 
diagonalizable and hence leave invariant some
complement subspaces $K^a_L$ and $K^b_L$ 
of $\ker \pi^{l+1,l}$ in $\bj^{l+1}(E)$. 
These complement subspaces do not depend on 
$L\in \R^*_+$ since $a_L$ and $a_{L'}$ 
(resp. $b_L$ and $b_{L'}$) commute 
for all $L, $ $L'\in \R^*_+$. 
We deduce that in a diagonalization basis 
of $a_L$, where the eigenvalues
are ordered in the decreasing way, $L^{-p}, L^{-p-q}, L^{-p-2q}, \cdots, 
L^{-p-(l+1)q}$, there
exists a lower unipotent endomorphism
$T$ such that $b_L = T\circ a_L \circ T^{-1}.$
It follows that 
$a_L\circ b_L^{-1} = (a_L \circ T \circ a_L^{-1})\circ T^{-1}$ is a product of unipotent endomorphisms 
$(a_L \circ T \circ a_L^{-1})$ and $T^{-1}$. The 
coefficients of $a_L\circ T\circ a_L^{-1}$ converge outside the diagonal 
to $0$
and the same holds for those of the product
$a_L\circ T\circ a_L^{-1}\circ T^{-1}.$
\epr
Note that  every one-parameter 
$(p,q)$-group of endomorphisms $(a_L)_{L\in \R^*_+}$
of jet bundle $\bj^l(E)$, $l\in \Nn$,
induces a one-parameter group of endomorphisms
of the bundle $\bj^l(E_{|\bh})$ denoted by $(a_L)_{L\in \R^*_+}$
too.

 Now, let $(U_L)_{L\in \Rp}$ be
a family of finite dimensional subspaces of $\gme$
which are asymptotically ample, meaning ample
for $L$ large enough. 
We equip them with scalar products $\cg \ , \ \cd_{L\in \R^*_+}$ 
For $L$ large enough the latter  induces after push-forward by
$gr^0$ and  $gr^1$ respectively, 
a sequence of Riemannian metrics
 $g_L^0$, $g_L^1$ on $E$ and $E\oplus (H^*\otimes E)$ respectively, 
 see Definition \ref{defpb}.  
It also induces the sequence of push-forwarded measures $j_{0*}\mu_{\bi^0}  $ 
 and $j_* \mu_{\bi^1_i} $
on $(H^\circ \oplus Sym^2 (H))\otimes E$.
\begin{defin}\label{tamed}
The family  $(U_L, \cg \ , \cd_L)_{L\in \R^*_+}$ is 
said to be $(p,q)$-tamed if and only if 
there exists a one-parameter $(p,q)$-group
of endomorphisms  $(a_L)_{L\in \Rp}$  of jet bundles such that
\begin{itemize}
\item When $n\geq 2$, 
$(a_L)^{-1*} g_L^1$
converges to a metric $g_\infty$ on 
$E\oplus (H^*\otimes E)$
and for every $\izn$, $(a_L)_* j_*\mu_{\bi_i^1}$
converges to a measure $\mu_\infty^i$. 
\item When $n=1$, 
$(a_L^*)^{-1} g^0_L$
converges to a metric $g_\infty$ on 
$E$
and  $(a_L)_* j_{0*}\mu_{\bi^0}$
converges to a measure $\mu_\infty$. 
\end{itemize}
\end{defin}
% % % % %
\begin{cor}\label{Corollaire 2}
Let $E$ be a real line bundle over a smooth
manifold equipped with a Morse function. 
Let $(U_L, \cg \ , \ \cd_L)_{L\in \R^*_+}$ be a family of 
asymptotically ample finite dimensional linear subspaces of $\gme$,  
which are $(p,q)$-tamed for some $p,q>0$. Then,
for every $\izn$, 
$$\frac{1}{L^{qn}}\mathbb E(\nu_i) \underset{L\to \infty}{\to}
 \frac{1}{\sqrt \pi^n}
\iint_{(H^\circ\oplus Sym^2_i (H))\otimes E}
|(\alpha, \beta)^* dvol_{g_\infty}| 
d\mu^i_\infty (\alpha, \beta)$$
weakly on $M$ when $n\geq 2$.
When $n=1$, 
$\frac{1}{L^{q}}\mathbb E(\nu) \underset{L\to \infty}{\to}
 \frac{1}{\sqrt \pi}
\int_{T^*M\otimes E}
|\alpha^* dvol_{g_\infty}| 
d\mu_\infty (\alpha).$
\end{cor}
% % % % %
\bpr
From Theorem \ref{Theoreme 1}, 
for every $L\in \Rp$, 
$$ \mathbb E(\nu_i)= \frac{1}{\sqrt \pi^n} 
\iint_{(H^\circ \oplus Sym^2_i(H))\otimes E}
|(\alpha, \beta)^* dvol_{g^1}| j_*\mu_{\bi_i^1}(\alpha,\beta).
$$
Let $(a_L)_{L\in \Rp}$ be the one-parameter $(p,q)$-group
of endomorphisms of jet bundles
such that $(a^{-1}_L)^*g^1_L$ 
converges to $g_\infty$ as $L$ 
grows to infinity and $(a_L)_*j_*\mu_{\bi^1_i}$
converges to $\mu_\infty^i$. Then, 
$$ dvol_{a_L^{-1*}g^1_L} = a_L^{-1*} dvol_{g^1_L}
= L^{p+(n-1)(p+q)} dvol_{g^1_L},$$
so that 
$dvol_{a_L^{-1*}g^1_L} \equiL
L^{-p-(n-1)(p+q)} dvol_{g_\infty}.$
We perform the substitution
$ a_L \alpha = \tilde \alpha $ and $a_L \beta = \tilde \beta$,
so that
$$ \mathbb E(\nu_i)\equiL 
L^{-p-(n-1)(p+q)}L^{p+q+(n-1)(p+2q)}
\frac{1}{\sqrt \pi^n}
\iint_{(H^\circ\oplus 
Sym^2_i (H))\otimes E}
|(\tilde \alpha, \tilde \beta)^* dvol_{g_\infty}| 
d\mu^i_\infty (\tilde \alpha, \tilde \beta)$$
since $(a_L\circ j)_*\mu_{\bi^1_i}\underset{L\to \infty}{\to} \mu^i_\infty.$
The proof in the case $n=1$ is similar. 
\epr
% % % % % % % % % % % % % % % % % % % % % % % % % % % % % %
% % % % % % % % % % % % % % % % % % % % % % % % % % % % % %
\section{Random eigensections of a self-adjoint elliptic operator}
The aim of this section is to prove Theorem \ref{Theoreme 0.1} 
and Corollary \ref{Corollaire 0}, see \S \ref{222}.
We first recall in \S \ref{paragraphe 2}
the asymptotic estimates of the derivatives of the
spectral function along the diagonal, 
which are needed to get these results from Remark \ref{Remarque a}. 
A proof of these estimates is given in Appendix \ref{A1}
while several basic definitions 
on pseudo-differential operators 
are recalled in Appendix \ref{paragraphe II.1}.

% % % % % % % % % % % % % %
\subsection{Asymptotic derivatives
of the spectral function along the diagonal}\label{paragraphe 2}

Under the hypotheses of Theorem \ref{Theoreme 0.1}, 
we assume $P$ to be positive, see Remark \ref{positif}
and for every $L\in \R_+^*$, we denote 
by $ e_L \in \Gamma (M\times M, E\boxtimes E)
$
the  spectral function of $U_L$, so that 
$$ \forall s\in U_L, 
\forall x\in M, 
s(x)= \int_M h_E (e_L(x,y), s(y)) |dy|,$$
compare Definition \ref{Definition 2}.
In particular, if $(s_1, \cdots, s_{N_L}) $ denotes an orthonormal
basis of $U_L$, then for every 
$ x,y\in M, e_L (x,y)= \sum_{i=1}^{N_L} s_i(x)s_i(y).$
The metric $h_E$ induces an
isomorphism between the restriction of 
$E\boxtimes E$ to the diagonal of $M\times M$ 
and the trivial line bundle over $M$  and under this isomorphism,
for every
$x\in M, e_L(x,x)= \sum_{i=1}^{N_L} h_E(s_i(x), s_i(x))>0.$
The dimension $N_L$ of $U_L$ then reads
$N_L = \int_M e_L(y,y)|dy|.$
The asymptotic behaviour of the spectral function $e_L$
along the diagonal is given by 
 Theorem \ref{Theoreme 2},
due to Carleman \cite{Carleman}
when $m=2$ and to G{\aa}rding \cite{Garding}
in general.
\begin{thm}[\cite{Carleman}, \cite{Garding}]\label{Theoreme 2}
Let $P$ be an  elliptic  pseudo-differential operator 
of order $m>0$, which is self-adjoint and  bounded from below, acting
on a real Riemannian line bundle  over
a smooth closed manifold $(M, |dy|)$
of positive dimension $n$. Let $\syp$ be the principal 
symbol of $P$ and $e_L$ be its spectral function, $L\in \R_+$. 
Then, for every $x\in M$, 
$$ e_L (x,x)\equiL
\frac{1}{(2\pi)^n} \int_{K_L}  |d\xi|,$$
where $|d\xi|$ denotes the measure on $T_x^*M$ induced 
by $|dy|$ and 
\begin{equation}\label{KL}
K_L = \{\xi \in T^*M \, | \, \si_P(\xi)\leq L\}.
\end{equation}
\end{thm}
Note that $K_1= K$, see (\ref{KP}).
In particular, the asymptotic given by Theorem \ref{Theoreme 2}
neither depends on the Riemannian metric of $E$, 
nor on the global geometry of $M$, it only depends on the 
measure $|dy|$ of $M$ at $x$ and on the symbol 
of $P$. 
\begin{rem}\label{Weyl}
Recall that Theorem \ref{Theoreme 2} recovers Weyl's theorem, which
computes the dimension 
$$\frac{1}{ L^{\frac{n}{m}}}N_L \toL \int_M c_0(y) |dy|,$$
see (\ref{c0}).
For example, when  $P$ stands for the Laplace-Beltrami operator
associated to some Riemannian metric on $M$,
this formula reads
\beq
\frac{1}{\sqrt L^{n}}N_L \toL	
 \frac{1}{(2\pi)^n}  Vol(\mathbb B_n) Vol_{g} M,
\eeq
where $Vol(\mathbb B_n)$ denotes the volume of the unit ball in $\R^n$,
see \S \ref{Laplace}.
\end{rem}
In order 
to apply the results of \S \ref{paragraphe 1},
we have to know in addition the asymptotic
of the partial derivatives of the spectral function $e_L$ 
along the diagonal.
This is the object of Theorem \ref{Theoreme 3}.
\begin{thm}\label{Theoreme 3}
Under the hypotheses of Theorem \ref{Theoreme 2},
let $Q_1$ and 
$Q_2$  be two differential operators on $E$ 
with principal symbols $\sigma_{Q_1}$ and $\sigma_{Q_2}$, of 
order $|\sigma_{Q_1}|$ and $|\sigma_{Q_2}|$, acting
on the first and second variables
of $e_L$ respectively. Then, for every $x\in M$,
\begin{equation} \label{Th3}
 Q_1Q_2 e_{L|(x,x)}
= 
\frac{1}{(2\pi)^n} \int_{K_L}
\sigma_{Q_1}(i\xi) \overline{\sigma_{Q_2}(i\xi)} |d\xi| +
O(L^{\frac{n+|\sigma_{Q_1}|+|\sigma_{Q_2}| - 1}{m}}),
\end{equation}
see (\ref{KL}).
\end{thm}
Theorem \ref{Theoreme 3}  is proved by L. H\"ormander
in \cite{Hormander}
when $Q_1$ and $Q_2$ are trivial, providing
the order of the error term in Theorem \ref{Theoreme 2}. 
It is written in \cite{SaVa} when $Q_1$ and $Q_2$ are of the same order,
see Theorem 1.8.5 of \cite{SaVa}, but we did not find
a reference for the general case, which we need here. 
In the particular case where $P$ is the Laplace-Beltrami operator, 
Theorem \ref{Theoreme 3} is proved in \cite{Bin}, see also \cite{Nicolaescu}.
We give in Appendix \ref{A1} a  proof of Theorem \ref{Theoreme 3} which follows
closely  \cite{Hormander}. 
Note that when $|\sigma_{Q_1}|$ and $|\sigma_{Q_2}|$ are not of the same parity, 
the main term of the right-hand side of (\ref{Th3}) vanishes
since for every $\xi \in T^*M, \ \sigma_P(-\xi)= \si_P(\xi)$ 
while the principal symbols $\sigma_{Q_1}$ and $\sigma_{Q_2}$
are homogeneous.
When $|\sigma_{Q_1}|= |\sigma_{Q_2}| \, mod (2)$, (\ref{Th3}) reads
$$ Q_1 Q_2 e_{L|(x,x)}\equiL
\frac{1}{(2\pi)^n}(-1)^{\frac{|\sigma_{Q_1}|-|\sigma_{Q_2}|}{2}} \int_{K_L}
\sigma_{Q_1}(\xi) \overline{\sigma_{Q_2}(\xi)} |d\xi|.
$$

% % % % % % % % % % % % % % % % % % % % % % % % % %

\subsection{Metrics on symmetric tensor algebras}\label{tensoral}
Let $V$ be a real vector space and $V^*$ be
its dual. 
For every $k\in \Nn$, we  denote by $Sym^k (V)$
the space of symmetric $k$-linear forms on $V^*$.
For every $q\in Sym^k(V)$ and every $\xi \in V^*$, 
we set  $q(\xi )= q(\xi, \cdots, \xi)$
and $q(i\xi)= i^k q(\xi)$. 
For every $l\in \Nn$, we set 
\beq 
S^l(V) &=& \bigoplus_{0\leq k\leq l}
Sym^k(V),\\
S^l_+ (V) &=& \{q\in S^l(V) \ | \ q(\xi) = q(-\xi) \},\\
S^l_- (V )&=& \{ q\in S^l(V) \ | \ q(\xi )= -q(-\xi)  \}.
\eeq 

\begin{lem}\label{Lemme 8}
Let $V$ be a real vector space and $l\in \Nn$. 
Let $K\subset V^*$ 
and $\mu$ be a positive finite measure on $K$ such that 
\begin{enumerate}
\item $-id$ preserves $K$ and $\mu$ 
\item \label{(2)}The support of $\mu$ is not included in 
any  degree $l$ algebraic hypersurface of  $V$.
\end{enumerate}
Then, the bilinear form  
\beq \kappa^l : 
S^l(V) \times S^l(V) & \to & \C \\
(q_1, q_2) & \mapsto & \frac{1}{\mu(K)} \int_K 
q_1(i\xi) \overline{q_2(i\xi) }d\mu (\xi)\in \C 
\eeq
 associated to $(K,\mu)$  only takes real values and defines a scalar product
on $S^l(V)$. Moreover, 
$S^l_+ (V)$ and $S^l_-(V)$ are orthogonal to each other with respect to $\kappa^l$.
\end{lem}
\bpr 
The form $\kappa^l$ is bilinear 
and the change of variables 
$\xi \in K \mapsto -\xi \in K $  yields that  
$S^l_+ (V)$ and $S^l_-(V)$ 
are orthogonal to each other. Moreover, the restrictions of $\kappa^l$ to 
$S^l_+(V)$ and $S^l_-(V)$ 
are real and symmetric, so that $\kappa^l$ itself 
is symmetric and takes only real values. Lastly,
if $q= \sum_{j=0}^{\lfloor l/2\rfloor }
q_j\in S^l_+(V),$
where for every $ j\in \{0, \cdots, \lfloor l/2\rfloor\},$
$ q_j \in Sym^{2j}(V^l)$, then 
$$\kappa^l(q,q)= \frac{1}{\mu(K)}
\int_K (\sum_{j=0}^{\lfloor l/2 \rfloor} (-1)^j q_j(\xi))^2 d\mu(\xi),$$
so that the restriction of $\kappa^l$ to $S^l_+(V)$
is non negative and the second hypothesis 
implies that it is positive definite. The same conclusion
holds for  the
restriction of $\kappa^l$ to $S^l_-(V)$, hence the result.
\epr
\begin{rem}\label{Remarque 1}
\begin{enumerate}
\item Under the hypotheses of Lemma \ref{Lemme 8}, 
the restriction of $\kappa^1$ to $Sym^1(V)= V$
defines a scalar product on $V$.
\item
If the measure $\mu$ can be chosen to be 
the absolute value of an alternated   $\dim V$-linear 
form on $V$, then the scalar products $\kappa^l$ given by 
Lemma \ref{Lemme 8} do not depend on the choice of this
form and only depend on $K$. This is the case when 
$K$ is bounded and has a non-empty interior. 
\end{enumerate}
\end{rem}

\subsection{Proof of Theorem \ref{Theoreme 0.1} and Corollary \ref{Corollaire 0}}\label{paragraphe 3'}

\subsubsection{Induced metric on the symmetric tensor bundle}\label{para 3.3.1}

Since $P$ is real and self-adjoint, the set
$K_L=\{\xi \in T^*M \, | \, \si_P(\xi)\leq L\}$
is invariant under $-Id$ and induces
thus a Riemannian metric on $M$
and
even 
on all symmetric tensor powers $S^l(TM)$, 
$l\in \Nn$, see Lemma \ref{Lemme 8} and Remark \ref{Remarque 1}.

\begin{defin}\label{glL}
For every $L\in \R^*_+$ and $l\in \Nn$, we
denote by $\kappa^l_L$ the Riemannian metrics induced by $K_L$
on $S^l(TM)$, see Lemma  \ref{Lemme 8}. 
\end{defin}
Together with the metric $h_E$, 
$\kappa^l_L$ induces a metric 
on $S^l(TM)\otimes E^*$
and by duality a metric on $S^l (T^*M)\otimes E$,
still denoted by $\kappa^l_L$. 
\begin{Proposition}\label{Proposition 1}
Under the hypotheses of Theorem \ref{Theoreme 2},
for every $l\in \Nn$
and every large enough $L\in \R^*_+$,
 $(U_L, \cg \, , \, \cd_L)$
is $l$-ample 
and $(\frac{n}{2m}, \frac{1}{m})$-tamed. Moreover,
the push-forward of $\cg \, , \, \cd_L$
under $j^l : \us_L \to \jle$ satisfies
$$ j^l_*\cg \, , \, \cd_L\equiL
L^{\frac{n}{m}} c_0 \kappa^l _L,$$
see \S  \ref{param}. 
\end{Proposition}
\bpr
From Lemma \ref{Lemme 4}, 
the push-forward $h_L$ of 
$\cg \, , \, \cd_L$ under $j^l$ induces on $\jle^*$ 
the metric 
$j^l j^l e_L$.
Let us fix a torsion-free  connection $\nabla$ on $TM$ and
a connection $\nabla^E$ on $E$. They induce a 
decomposition $\bj^l (E)\cong 
%\oplus^{l}_{k=0}
S^l (T^*M )\otimes E$
which 
equips $\bj^l(E)^*$ with 
the metric $\kappa^l_L$. From Theorem \ref{Theoreme 3} follows that 
the metrics $h^1_L$ and $L^{\frac{n}{m}}c_0 \kappa^l_L$
are equivalent as $L$ grows to infinity. In particular the asymptotic value of 
the induced metric 
$ L^{\frac{n}{m}}c_0 \kappa_L^l$ on ${\cal J}^l (E)^*$ does not depend on the chosen decomposition $\bj^l (E)\cong 
S^l (T^*M )\otimes E$, 
see Lemma \ref{Lemme 6}. Now,
 $\kappa^l_L$ is $(p,q)$-tamed
with $p=n/(2m)$ and $q=1/m$.
Indeed, the one-parameter $(p,q)$-group 
of fibre bundles endomorphisms 
\beq a_L : 
\bigoplus^{l}_{k=0}
Sym^k (T^*M)\otimes E &\to & \bigoplus^{l}_{k=0}
Sym^k (T^*M)\otimes E\\
(q_k)_{k\in \{0, \cdots , l\}} & \mapsto & (L^{-\frac{n}{2m}- \frac{k}{m}} q_k)_{k\in \{0, \cdots , l\}}.
\eeq
is such that $L^{n/m}a_L^{-1*}\kappa_L^l$ converges to the metric
associated to $ (K, d\xi)$ 
given by Lemma \ref{Lemme 8}.
\epr
\begin{cor}\label{Corollaire 3}
Under the hypotheses of Theorem \ref{Theoreme 0.1},
the push-forward of $\cg \ , \ \cd_L$ 
under $j$ gets equivalent, as 
$L$ grows to infinity and when $n\geq 2$, to 
\beq
\big((H^\perp \times Sym^2 (H^*))\otimes E^*\big)^2
& \to & \R \\
((a_1, b_1), (a_2, b_2)) & \mapsto & 
\frac{1}{(2\pi)^n}
\big(\int_{K_L} h_E(a_1(\xi),a_2(\xi))
+ h_E(b_1(\xi) , b_2(\xi)) |d\xi| - \cdots \\
&& \cdots \frac{1}{\int_{K_L} |d\xi|}
\iint_{K_L^2} h_E(b_1(\xi),b_2(\xi')) |d\xi|
|d\xi'|
\big).
\eeq
When $n=1$, the push-forward of $\cg \ , \ \cd_L$ 
under $j_0$ gets equivalent, as 
$L$ grows to infinity, 
to 
$
(a_1, a_2)\in (T^*M\otimes E^*)^2  \mapsto 
\frac{1}{2\pi}
\int_{K_L} h_E(a_1(\xi),a_2(\xi))|d\xi|.
$
\end{cor}
In Corollary \ref{Corollaire 3}, 
%$K_L = \{\xi \in T^*M \ | \ \si_P (\xi)\leq L\}$
%and 
$H^\perp$
denotes the orthogonal of $H$ with respect
to the Riemannian
metric of $M$ associated to $K_L$,
given by Definition \ref{glL}. The 
distribution $H$ is defined in \S \ref{paragraphe 1}
and $j$ in \S \ref{paragraphe 1.1}. 
\bpr
From Proposition \ref{Proposition 1}, 
the metric $j^2 \# (j^2)^*$
of $\bj^2(E)^*$ 
gets equivalent to $L^{\frac{n}{m}} c_0 \kappa^2_L$
as $L$ grows to infinity. 
By restriction to the fibre product
$ \big(\bj^1(E)\times _{\bj^1(E_{|\mathcal H})} \bj^2(E_{|\mathcal H})\big)^*$,
we deduce that the metric induced on this space
gets equivalent to 
\beq
((\R \oplus TM \oplus Sym^2(H))\otimes E^*)^2
&\to &\R\\
((c_1, a_1, b_1), (c_2, a_2, b_2)) & \mapsto & 
\frac{1}{(2\pi)^n} 
\int_{K_L} h_E(c_1 ,c_2)(\xi)  - h_E(c_1 ,b_2)(\xi) - \cdots  \\
&&\cdots h_E( b_1,c_2)(\xi)
+ h_E(b_1, b_2 )(\xi)+ h_E(a_1, a_2)(\xi) |d\xi|
.
\eeq
We apply then Lemma \ref{Lemme 5} and Remark \ref{Remarque a}
to $F = \us_L$, $G= (\bj^1(E)\times_{\bj^1(E_{|\mathcal H})}\bj^2(E_{|\mathcal H}))^*$, $K_F = \bi^1$ 
and 
$K_G = (H^\perp \oplus Sym^2(H) )\otimes E^*$,
where the middle term $TM$ splits as $H\oplus H^\perp$.  
We deduce that the factors 
$H^\perp \otimes E^*$ and 
$Sym^2(H)\otimes E^*$ get asymptotically orthogonal, 
that the metric induced on $H^\perp \otimes E^*$ 
is asymptotically equivalent to $\frac{\mu(K_L)}{(2\pi)^n } $ times the one induced by $K_L$
and finally that the one induced on $Sym^2(H)\otimes E^*$
is equivalent to 
$$ (b_1, b_2) \mapsto \frac{1}{(2\pi)^n} \big( \int_{K_L} h_E(b_1, b_2)(\xi) |d\xi| - \frac{1}{\int_{K_L} |d\xi|}
\iint_{K_L\times K_L} h_E(b_1(\xi),b_2(\xi') ) |d\xi| |d\xi'|\big).$$
Indeed, with the notations of Lemma \ref{Lemme 5}, 
$L_G = (\R \oplus H )\otimes E^*$ gets a metric $c\#c^*$
for which the factors $\R \otimes E^*$ and $H\otimes E^*$ are
asymptotically orthogonal to each other and the metric on 
$\R \otimes E^*$ is $\frac{1}{(2\pi)^n} \mu(K_L) h_E$. Moreover,
the correlation $b\# c^*$ only involves the factors $\R \otimes E^*$
and $Sym^2 (H)\otimes E^*$ and reads
$$(c_1,b_2)\in E^*\oplus (Sym^2(H)\otimes E^*) \mapsto - \frac{1}{(2\pi)^n} \int_{K_L} h_E (c_1, b_2)(\xi) |d\xi|.$$
Finally $a\# a^*+ b\# b^*$ is a metric on $(H^\perp \oplus Sym^2(H))\otimes E^*$ 
for which both factors are asymptotically orthogonal, 
the metric induced on $H^\perp \otimes E^*$ is asymptotically 
equivalent to $\frac{\mu(K_L)}{(2\pi)^n }$ times the one induced by $K_L$,
and the one induced on $Sym^2(H)\otimes E^*$ is 
$$ (b_1, b_2) \in Sym^2(H)\otimes E^* \mapsto \frac{1}{(2\pi)^n} \int_{K_L} h_E(b_1,b_2)(\xi)|d\xi|.$$
We deduce now that
the correlation term $b\# c^* (c\# c^*)^{-1}c \# b^*$
just reads 
$$\frac{1}{(2\pi)^n \int_{K_L}|d\xi|} \iint_{K_L\times K_L} h_E(b_1(\xi), b_2(\xi')) |d\xi| |d\xi'|.$$
Hence the result.
\epr
\subsubsection{Proof of Theorem \ref{Theoreme 0.1}
and Corollary \ref{Corollaire 0}}\label{222}

We know from Proposition \ref{Proposition 1} that
$U_L = \bigoplus_{\lambda \leq L} \ker (P-\lambda Id)$
equipped with the $L^2$-scalar product 
$\cg \, , \, \cd_L$ gets ample for $L$ large enough
and $(\frac{n}{2m}, \frac{1}{m})$-tamed, see 
Definition \ref{tamed}.
From Corollary \ref{Corollaire 2}, we deduce that
$\frac{1}{L^{\frac{n}{m}}}\mathbb E (\nu_i)$ weakly converges 
on the whole  $M$
to the measure 
\begin{equation}\label{mes}
\frac{1}{\sqrt \pi^n}
\iint_{(H^\perp\times Sym^2_i (H^*))\otimes E}
|(\alpha, \beta)^* dvol_{g_\infty}| d\mu_\infty^i (\alpha, \beta),
\end{equation}
where the metric $g_\infty$ and the measure $\mu_\infty^i$
are given by Definition \ref{tamed}.
From Proposition \ref{Proposition 1} and Corollary 
\ref{Corollaire 3},
the factors $E$ and $H^*\otimes E$ are
orthogonal to each other with respect to $g_\infty$,
and $g_\infty$ restricts to $c_0 h_E$ on $E$ and
to the metric $g_P \otimes h_E$ on $H^*\otimes E$,
see (\ref{gp}). Likewise, from Corollary \ref{Corollaire 3}
the measure $\mu_\infty^i$ is a product of
the measure on $H^\circ \otimes E$ induced by $g_P$
and $h_E$,
and the measure on $Sym^2_i (H)\otimes E$ induced
by (\ref{sym}) and $h_E$. We deduce that 
$dvol_{g_\infty}= \frac{1}{\sqrt c_0} dvol_{h_E}$
and that (\ref{mes}) becomes 
$$\frac{1}{\sqrt \pi^n \sqrt c_0} \mathbb E(i,\ker dp) 
\Big(\int_{H^\perp\otimes E} |\alpha| d\mu_{P}(\alpha)\Big)|dvol_{P}|.$$
We conclude thanks to the equality 
$$ \int_{H^\perp\otimes E} |\alpha| d\mu_{P}(\alpha)
 = \int_\R |a| e^{-a^2} \frac{da}{\sqrt \pi}= \frac{1}{\sqrt \pi}.$$
 When $n=1$, 
 $\frac{1}{L^{\frac{1}{m}}}\mathbb E (\nu)$ weakly converges 
 to the measure 
 \beq \frac{1}{\sqrt \pi}
 \int_{T^*M\otimes E}
 |\alpha^* dvol_{g^0_\infty}| d\mu_\infty (\alpha)&=& \frac{1}{\sqrt \pi \sqrt c_0}
 \int_{T^*M\otimes E} |\alpha| d\mu_{K}(\alpha)|dvol_P|\\ 
 &=&\frac{1}{\pi \sqrt c_0} |dvol_{P}|. \square
 \eeq

\bpr[ of Corollary \ref{Corollaire 0}]
It  is a consequence of Theorem
\ref{Theoreme 0.1} after integration of the 
constant function 1, compare Corollary \ref{Corollaire 1}.
\epr
% % % % % % % % % % % % % % % % % % % % % % % % % % % % % % % % % % % 

% % % % % % % % % % % % % % % % % % % % % % % % % % % % % % % % % % % %
\section{Examples}\label{exemples}
We investigate in this third section two examples, 
the Laplace-Beltrami operator in 
\S \ref{Laplace}, where we prove 
Corollary \ref{Corollaire 0'} and Proposition \ref{prop}, and
the Dirichlet-to-Neumann operator in \S \ref{Dirichlet},
where we prove Corollary \ref{d2n}.

\subsection{The Laplace-Beltrami operator}\label{Laplace}
\subsubsection{Proof  of Corollary \ref{Corollaire 0'}}
The principal symbol of 
the Laplace-Beltrami operator $\Delta_g$ reads 
$ \sigma_{\Delta_g}: \xi \in T^*M \mapsto g(\xi,\xi) \in \R,$
so that the  compact $K$ defined by (\ref{KP}) reads 
$$ K = \{\xi \in T^*M\, | \, g(\xi,\xi)\leq 1\}.$$
The Riemannian metric $g_{\Delta_g}$ induced on $M$ by the pair $(K, |d\xi|)$
reads at every point $x\in M$,
$ (u,v)\in T_xM^2 \mapsto \frac{1}{(2\pi)^n} \int_{K}
\xi (u)\xi (v) |d\xi|$ by (\ref{gp}),
so that 
\begin{equation}\label{gdel}
 g_{\Delta g}= c_1 g 
 \end{equation}
 and
 \begin{equation}\label{vol}
 |dvol_{\Delta_g}|= \sqrt {c_1}^n |d\xi|,
 \end{equation}
 where 
 \begin{equation}\label{c1}
c_1 = \frac{1}{(2\pi)^n } \int_{K}
  \xi_1^2 |d\xi|. 
  \end{equation}
Let us choose an  orthonormal basis
 $(\partial/\partial_{x_1}, \cdots, \partial/\partial_{x_n})$
 of $T_xM$ such that 
  $(\partial/\partial_{x_1}, \cdots, \partial/\partial_{x_{n-1}})$
spans $H_x$ and let us denote by 
  $(\xi_1, \cdots, \xi_n)$ its dual basis.
%Together with an orthonormal basis of $E_x$, 
They induce  isomorphisms
 $ Sym^2(H)\cong Sym(n-1,\R)$
 and $Sym^2(H)^*  \cong Sym(n-1,\R)^*$.
 From Corollary \ref{Corollaire 3}, when $n>2$ 
 the metric induced by $(K, |d\xi|)$
 on $Sym(n-1, \R)^*$
then  reads 
$$
\forall (A,B) = ((a_{ij})_{1\leq i,j\leq n-1},(b_{ij})_{1\leq i,j\leq n-1}) \in (Sym^*(n-1,\R))^2,
$$
$$
\cg A, B\cd_{{\Delta_g}} =
\upn
\big(\int_{K} A(\xi)B(\xi) |d\xi|  - 	
\frac{1}{\int_{K} |d\xi|}
\int_{K}  A(\xi) |d\xi|  
\int_{K} B(\xi) |d\xi|  \big)
$$
where 
$$\int_{K} A(\xi)B(\xi) |d\xi| = 
 \int_{K}
(\sum_{i=1}^{n-1} a_{ii} \xi_i^2 + 2\sum_{1\leq i<j \leq n-1} a_{ij} \xi_i \xi_j)
(\sum_{i=1}^{n-1} b_{ii} \xi_i^2 + 2\sum_{1\leq i<j \leq n-1} b_{ij} \xi_i \xi_j)
|d\xi| 
$$
and 
\beq
\int_{K}  A(\xi) |d\xi|  
\int_{K} B(\xi) |d\xi|   &=& 
\int_{K}(\sum_{i=1}^{n-1} a_{ii} \xi_i^2 + 2\sum_{1\leq i<j \leq n-1} a_{ij} \xi_i \xi_j)|d\xi|\cdots\\
&&\cdots \int_{K}(\sum_{i=1}^{n-1} b_{ii} \xi_i^2 + 2\sum_{1\leq i<j \leq n-1} b_{ij} \xi_i \xi_j )|d\xi|,
\eeq
so that 
\beq
\cg A, B\cd_{K}
 & = &(c_4 -\frac{c^2_1}{c_0})\sum_{i=1}^{n-1} a_{ii}b_{ii}
 + (c_2 - \frac{c^2_1}{c_0}) \sum_{1\leq i\not=j\leq n-1} a_{ii}b_{jj}+ 4c_2 \sum_{1\leq i<j\leq n-1} a_{ij}b_{ij}\\
& = & 2c_2 \Tr (AB)+ (c_2 - \frac{c^2_1}{c_0}) (\Tr A)(\Tr B) ,
\eeq
where
$$\begin{array}{lcl}
c_4&=& \upn \int_{K} \xi_1^4  |d\xi|,\\
c_2 &=& \upn \int_{K} \xi_1^2 \xi_2^2 |d\xi| \text{ and } \\
c_0 &=& \upn \int_{K}  |d\xi| .
\end{array}
$$
This indeed follows from 
the relation $c_4 = 3c_2$,
see \cite{Bin}, \cite{Nicolaescu} and from the fact
that $\int_{K} \xi_1^k \xi_2^l |d\xi| = 0$ 
whenever $k $ or $l$ is odd. 
Note that
\begin{equation}\label{c2}
 \begin{array}{lcl}
c_2 &=&\frac{c_0}{(n+4)(n+2)},\\
c_1 &= &\frac{c_0}{n+2} \text{ and }\\
c_2 - \frac{c^2_1}{c_0} &=& \frac{-2c_2}{n+2},
\end{array} 
\end{equation}
  see \cite{Bin} and \cite{Nicolaescu}.
Hence, the scalar product induced by $(K, |d\xi|)$ 
on $Sym(n-1, \R)^*$ is given, with the notations 
of the appendix B of \cite{Nicolaescu}, by the symmetric endomorphism
$ 2c_2 Q(a,b,c)$ with $a= \frac{n+1}{n+2}$, $b= \frac{-1}{n+2}$
and $c=1$. As a consequence, 
 the induced  scalar product on $Sym(n-1, \R)$
is given by the symmetric endomorphism $\frac{1}{2c_2} Q(a',b',c')$
with $a'= \frac{4}{3}$, $b'= \frac{1}{3}$ and $c'= 1$,
see \cite{Nicolaescu}.
Hence, for every $(A,B)\in 
Sym(n-1,\R)^2$,
$$\cg A,B \cd_{\Delta_g} =\frac{1}{2c_2} (\Tr (AB)+ \frac{1}{3}(\Tr A ) (\Tr B)).
$$
Finally,
\beq
\mathbb E(i, \ker dp) & = & \int_{Sym^2(H)}|\det \beta | 
d\mu_{\Delta_g} (\beta) \\ 
& = & \frac{1}{c_1^{n-1}} \int_{Sym(i,n-1-i,\R )} |\det B|
e^{- \frac{1}{2c_2}(\Tr (B^2) + \frac{1}{3} (\Tr B)^2 )} d\mu_{\Delta_g}(B) \\ 
& = & \frac{\sqrt {c_2}^{n-1}}{c_1^{n-1}} \mathbb E(i,n-1-i),
\eeq
see (\ref{ein}), since from (\ref{gdel}), $|\det B| = c_1^{n-1}|\det \beta|
$ under the substitution $B = \beta$.
We deduce from Theorem \ref{Theoreme 0.1} 
and (\ref{vol})
the weak convergence on $M$ 
$$ \frac{1}{\sqrt L^n } \mathbb E(\nu_i) \underset{L\to \infty}{\to}
 \frac{1}{ \sqrt \pi^{n+1} \sqrt{c_0}}\frac { \sqrt {c_2}^{n-1}}
{c_1^{n-1}} \mathbb E(i,n-1-i)\sqrt {c_1^n} |dvol_g|.$$
The result follows now from (\ref{c2}) and Corollary \ref{Corollaire 0}. 
The proof goes along the same lines when $1\leq n\leq 2$
and the result remains true in these cases.\ $\square$

% % % % % % % % % % % % % % % % % % % %
% % % % % % % % % %
\begin{exa}
When $n= 2$, $\mathbb E(0,1)= \mathbb E(1,0)= \int_0^{+\infty} a e^{-\frac{2}{3} a^2} d\mu(a)=
\frac{\sqrt 3}{2\sqrt 2\sqrt \pi}$, so that
from Corollary \ref{Corollaire 0'}, for every $j\in \{0,1\}$,
\begin{eqnarray}\label{b0m}
\frac{1}{L} \mathbb E(\nu_j) &\toL &\frac{1}{8\pi^2}|dvol_g|\\
\text{ and }\limsup_{L\to \infty} \frac{1}{L} \mathbb E(m_j)&\leq &\frac{1}{8\pi^2} Vol_g(M).
\end{eqnarray}
\end{exa}
% % % % % % % % % % % % % % % % % % % % % % % % % % % % % % %
\subsubsection{Proof of Proposition \ref{prop} }\label{Aprop}
By Corollary \ref{Corollaire 0'} and Weyl's Theorem, see Remark \ref{Weyl}, 
it is enough to prove that there exist $C>0$ and $\delta >0$ such that
$$\forall n\in \Nn,\  \sum_{|\frac{i}{n}-\frac{1}{2}|\geq \epsilon} \mathbb E(i,n-i)\leq C\exp (-\delta n^2),$$
since $\log Vol (\mathbb B^n )\sim_{n\to \infty}- \frac{n}{2}\log n$. 
Now, if 
$ d\mu_{GOE}$ denotes
the Gaussian probability measure on $Sym(n,\R)$ associated to 
the scalar product $\cg A,B\cd = Tr (AB)$, then
the Gaussian probability measure $\mu$
associated to (\ref{scal3}) satisfies the bound
$\mu \leq c_n\mu_{GOE}$
with $c_n = O(n).$
Indeed, $\frac{1}{2}\Tr A^2 + \frac{1}{6}(Tr A)^2 \geq  \frac{1}{2}\Tr A^2$,
whereas the ratio between the determinants of these scalar product
is a $O(n)$, see (B.6) in \cite{Nicolaescu}.
Now, Theorem 1.6 of \cite{GaWe3} provides the result.

\subsection{The Dirichlet-to-Neumann operator}\label{Dirichlet}
Let $(W,g)$ be a smooth compact Riemannian manifold 
with boundary
and $\Delta_g$ be its Laplace-Beltrami operator.
Let us denote by $M$ the boundary of $W$ and
for every  smooth function $f :  M \to \R$,
we denote by $u\in C^\infty(M,\R)$ the solution of the Dirichlet
problem 
$$ \left\{
\begin{array}{lcl}
\Delta_g u &=& 0\\
u_{| M}&= & f.
\end{array}
\right.
$$
We then denote by $\partial_n u :  M \to \R$
the outward normal derivative of $u$ along $M$. 
Then, the Dirichlet-to-Neuman operator 
$\Lambda_g$ reads
\beq 
\Lambda_g : C^\infty (M, \R )&\to& C^\infty ( M, \R )\\
f & \mapsto& \partial_n u.
\eeq
\begin{thm}[\cite{dirichlet}]\label{dton}
Let $(W,g)$ be a smooth compact Riemannian manifold
with boundary $M$. 
The Dirichlet-to-Neumann operator $\Lambda_g$ is an elliptic pseudo-differential operator of order one on $M$. 
Its principal symbol equals  $\xi\in T^*M \mapsto  \|\xi\|_g.$
\end{thm}
\bpr[ of Corollary \ref{d2n}]
The compact $K_\Lambda$ defined by (\ref{KP})
coincides with  $K_{\Delta_g}$, where $K_{\Delta_g}$ is
the compact associated to the Laplace-Beltrami
operator on $M$ induced by the restriction
of $g$ to $M$. The proof of Corollary \ref{d2n} thus
goes along the same lines as the one of Corollary \ref{Corollaire 0'}.
\epr

% % % % % % % % % % % % % % % % % % % % % % % % % % % % % % % % % % % % % % % % % % % % % % % % % % % % % % % % %

\section{Some related problems}\label{res}
Let us mention several related problems which we plan
 to discuss in a separate paper.
First, we may consider, as our probability space, the 
span of eigensections with  eigenvalues belonging to a window 
$[a(L)L, L]$ instead of $[0,L]$, where $a$ is some function of $L$, 
compare \cite{LerarioLundberg},
\cite{SarnakWigman}.
That is, we may set 
 $$U_L^a = \bigoplus_{\lambda \in [a(L)L, L]} \ker (P - \lambda Id).$$
When $\lim_{L\to \infty} {a(L)} = \gamma\in [0,1]$,
Theorem \ref{Theoreme 0.1} still holds true, with the following modifications:
$K$ given by (\ref{KP}) should be replaced 
by the annulus $K^\gamma= \{\xi \in T_x^*M \, | \, \gamma \leq \sigma_P (\xi)\leq 1\}$
and when $\gamma=1$, we should assume that 
$L^{-\frac{1}{m}} = o(1-a(L))$
and replace $|d\xi|$ by some Lebesgue measure on the sphere $K^1$. 
In the latter case for example, when $P$ stands for the
Laplace-Beltrami operator associated to some Riemannian metric
$g$ on the closed $n$-dimensional manifold $M$, we get the weak 
convergence
$$ \frac{1}{\sqrt L^n} \mathbb E(\nu_i) \tol 
 \frac{1}{\sqrt {\pi}^{n+1}}\frac {  1}
 {\sqrt {n ( n+2)^{n-1}}}\mathbb E_S(i,n-1-i) |dvol_g|,$$
 where
$
 \mathbb E_S(i,n-1-i)= \int_{Sym(i,n-1-i), \R} |\det A | d\mu_S(A),
 $
and $\mu_S$
  is the Gaussian measure on $Sym(n-1,\R)$ associated to
  the scalar product 
 \begin{equation}\label{scalS}
  (A,B)\in Sym(n-1,\R)^2 \mapsto \frac{1}{2}\Tr (AB)+ \frac{1}{2}(\Tr A)(\Tr B)\in \R,
  \end{equation}

 Finally, a manifold of special interest is the round unit sphere, where 
 we may consider the space of pure harmonics 
 $ U^1_L = \ker (P-L Id)$
 as a probability space, compare \cite{NazarovSodin}, \cite{LerarioLundberg}.
 Recall that the spectrum of the Laplace-Beltrami operator on the round
 unit $n$-dimensional sphere is the set $\{l(l+n-1)\, | \, l\in \Nn\}$ 
 and that the eigenspace associated to the eigenvalue $\lambda_l = l(l+n-1)$
 has dimension  ${n+l \choose n}- {n+l-2 \choose n}$.
This case of pure spherical harmonics is unfortunately not a special case 
of the previous one, because $\gamma = 1$
but $L^{-1/m}$ cannot be a $o(1-a(L))$. 
However, the result remains valid and we also get the weak convergence
$$ \frac{1}{\sqrt L^{n} } \mathbb E (\nu_i) \underset{l\to \infty}{\to} 
\frac{\mathbb E_S(i,n-1-i)}{\sqrt  \pi^{n+1}\sqrt{
n  (n+2)^{n-1}}}
  |dvol_g|$$
 on the whole $M$. In the case 
 $n= 2$, this provides the upper estimate
\begin{equation}
\limsup_{l\to \infty }\frac{1}{L} \mathbb E(b_0) \leq \frac{1}{\pi \sqrt 2 },
\end{equation}
for the expected number $b_0$ of connected component 
of pure spherical harmonics, 
compare relation (2.41) of \cite{Nicolaescu}.

% % % % % % % % % % % % % % % % % % % % % % % % % % % % % % % % % % % % % % % % % % % % % % % % % % % % % % % % %
% % % % % % % % % % % % % % % % % % % % % % % % % % % % % % % % % % % % % % % % % % % % % % % % % % % % % % % % %
\appendix

\section{Appendix }\label{AppB}

\subsection{The incidence varieties}\label{A0}
We recall that for every subspace $U$ of $\gme$, 
\beq 
\Delta_0 &=& \{s\in U\, | \ s \text{ does not vanish transversally} \} \text{ and}\\
\Delta_1 &=& \Delta_0 \cup \{s\in U\setminus \Delta_0\ | \ p_{|s^{-1}(0)} \text{ is not Morse}, \}
\eeq
see \S \ref{paragraphe 1.1} (\ref{Delta1}).
\begin{lem}[compare Proposition 2.8 of \cite{GaWe4}]\label{Lemme 1}
Let $E$ be a real line bundle 
over a smooth manifold $M$
equipped with a Morse function $p  : M\to \R$
and let $U$ be a relatively $l$-ample  linear subspace of $\gme$,
$l\in \{0,1\}$. Then, $\bi^l$ is a submanifold of 
$\us_{|M\setminus Crit(p)} = (M\setminus Crit(p))\times U$
of  codimension 
$rank (\jleh)$. 
Moreover,
$\Delta_0$ 
coincides with the critical locus of $\pi_U : \bi^0 \to U$,
whereas 
$\Delta_1\setminus \Delta_0$
coincides with the critical locus of 
the restriction 
 $\pi_{U|(\bi^1\setminus \pi_U^{-1}(\Delta_0))}
:\bi^1\setminus \pi_U^{-1}(\Delta_0) \to U$.
\end{lem}
From Lemma \ref{Lemme 1} and Sard's Lemma, 
when $U$ is relatively $l$-ample, $l\in\{0,1\}$, $\Delta_l$
has measure zero.
\bpr
Let us first assume that $l=0$ and let $(x,s)\in \bi^0$.
We fix some connection $\nabla^E$ on $E$. 
Then, the differential of $j^0$ at $(x,s)$
reads
\beq
d_{|(x,s)}j^0 : T_{(x,s)}\us \to T_{(x,0)} E\\
\xsd \mapsto (\dot{x}, \dot{s}(x)+ \nabla^E_{\xd} s ).
\eeq 
Since $j^0$ is onto, 
$d_{|\xs} j^0 $ is onto as well 
and it follows from the implicit function theorem
that $\bi^0$ is a codimension one  submanifold of $\us_{|M\setminus Crit(p)}$ 
with tangent space 
\bq\label{TI0}
 T_{\xs} \bi^0 = \{\xsd \in T_{\xs} \us \ | \ \sd(x)+ 
\nae_{\xd} s = 0  \}.
\eq
Moreover, the differential 
$ d_{|\xs} \pi_U : \xsd \in T_{\xs} \bi^0 \mapsto
\sd \in T_s U = U $
is onto if and only if  $\nae s $ is, since $j^0$ is onto. 
Hence, $\Delta_0$ coincides
with the locus of the singular values of $\pi_U : \bi^0 \to U$.

Now, assume that $l=1$ and let $\xs \in \bi^1$.
The differential of $j^1_\bh$ at $\xs$ reads
\beq
d_{|\xs} j^1_\bh : T_{\xs} \us &\to& T_{(x,0)} \bj^1 (E_{|\bh} )\\
\xsd &\mapsto  & (\xd, j^1_\bh (\sd)+ \nabla^\bj_{\xd} (j^1_\bh (s))),
\eeq
where $\nabla^\bj$ denotes a connection
on the bundle $\bj^1 (E_{|\bh} )$. Since 
$j^1_\bh$ is onto,  $d_{|\xs} j^1_\bh$ 
is onto as well and it follows from the implicit function theorem
that $\bi^1$ is a submanifold of $\us_{|M\setminus Crit(p)}$ 
of codimension $rank(\bj^1 (E_{|\bh} ))= n$, with
tangent space 
\bq\label{TI1}
T_{\xs} \bi^1 = \{\xsd \in T_{\xs} \us \ | \
j^1_\bh (\sd)+ \nabla^\bj_{\xd} (j_\bh^1(s))= 0  \}.
\eq
Let us assume that $s\notin \Delta_0$ and let $\xsd
\in \ker d_{|\xs} \pi_U.$
Then $\sd=0$, which implies that
$\nabla_{\xd}^E s= 0$, so that 
$\xd \in \ker \nabla s_{|x}= H_x$. 
Then, $ 0= \nabla_{\xd}^\bh (j_\bh^1 (s))=
j^2_\bh (\xd, \cdot )$,
so that $\xd \in \ker j^2_\bh (s).$
We deduce that the kernel of $d_{|\xs} \pi_U$
is reduced to  $\{0\} $ if and only if  $j^2_\bh$
is non-degenerate. 
From Lemma \ref{Lemme 2}, $j_\bh^2(s)$
is non-degenerate if and only if  $s\notin \Delta_1$.
\epr
\begin{rem}\label{rema}
It follows from the proof of Lemma \ref{Lemme 1}
that for every $s\in \bi^1\setminus \Delta_1$, the operator
$\nabla^{\bj}(j_\bh^1 (s))$ which appears in (\ref{TI1})
is invertible. 
\end{rem}
% % % % % % % %
\begin{lem}\label{Lemme 2}
(compare Lemma 2.9 of \cite{GaWe4})
Let $E$ be a real fibre bundle over a smooth
manifold $M$ equipped
with a Morse function $p : M\to \R$. Let $s$
be a section of $E$ which vanishes transversally 
and $x\in M\setminus Crit(p)$ be a critical point of $p_{|s^{-1}(0)}.$ 
Let $\lambda \in E^*_x$ such that $\lambda \circ \nae s_{|x}
= d_{|x}p.$ Then,
$$ \lambda \circ \nabla^p (\nae s_{|\bh_x})_{|x} = 
\lambda \circ \nabla (\nae s)_{|x} - \nabla(dp) = -
\nabla^s (dp_{|s^{-1}(0)}).$$
\end{lem}
% % % % % % % %
In Lemma \ref{Lemme 2}, $\nae$,  $\nabla^p$, $\nabla^s$
and  $\nabla$
denote connections on, respectively,
 the fibre bundles $E$,  $H$, $T(s^{-1}(0))$ and  $TM$. 
These connections induce  connections on, respectively,
 $H^*\otimes E$, 
$T^*(\suz)\otimes E$ and $T^*M\otimes E$,
denoted in the same way by $\nabla^p $, $\nabla^s $ 
and $\nabla$.
Note that $\nae s$, $\nabla^p (\nae s_{|\bh})$ 
and $\nabla^s (dp_{|s^{-1}(0)})_{|x}$ 
do 
not depend on the choices of $\nae$, $\nabla^p$, $\nabla^s$,
whereas $\nabla (\nae s) $ and $\nabla dp$ depend on the choice of $\nabla$. 
% % % 
\bpr
Let $v,w$ be two vector fields on $\suz$ defined in the neighbourhood of $x$.
Then, 
$$ 0 = \nae_v (\nae_w s)_{|x} = \nabla( \nae s)(v,w)+ \nae_{\nabla_v w} s$$
and likewise 
$ \nabla^s (dp)_{|x}(v,w) = d_{|x} (dp(w))(v)= \nabla (dp)
(v,w)+ d_{|x}p (\nabla_v w).$
We deduce the relation 
$\nabla_{|x} (dp_{|\suz})(v,w)= \nabla (dp)_{|x} (v,w) - \lambda \circ
\nabla (\nae s) (v,w).$
Likewise, if $v'$ and $w'$ are two 
vector fields of $\bh_x$ defined in the neighbourhood of $x$,
we have 
$$ 0 = d_{|x} (dp(w'))(v')= \nabla (dp)(v',w')+ dp(\nabla_{v'}w')
$$
and 
$ \nabla^p(\nae s)(v',w') = \nabla^E_{v'}(\nae_{w'} s)= \nabla(\nae s) (v',w')+
\nae_{\nabla_{v'}w'} s.$
Finally, 
$$ \lambda \circ \nabla^p(\nae s)_{|x} = \lambda \circ \nabla (\nae s)_{|x}
- \nabla (dp)_{|x} = - \nabla^s (dp_{|\suz}).$$
\epr

% % % % % % % % % % % % % % % % % % % % % % % %
% % % % % % % % % % % % % % % % % % % % % % % % % %

% % % % % % % % % % % % % % % % % % % % % % % % % % % %
\subsection{Pseudo-differential operators}\label{paragraphe II.1}

Let $M$ be a smooth manifold of positive dimension $n$ and $E $ be a real line bundle  over $M$.
We denote by $\Gamma(M,E)$ the space of smooth global sections of $E$.
\begin{defin}\label{Definition 4.0}
(compare Definition 18.1.32 of \cite{HormanderPDEIII})
%see \cite{Hormander-psd}, \cite{Hormander}, \cite{Hormander-Fourier})
A linear operator $P : \Gamma(M,E)\to \Gamma(M,E)$ 
is called pseudo-differential of order $m\in \R$ 
if and only if  there exist an atlas  $(U_i)_{i\in I}$ of $M$ 
and local trivializations $\Phi_i : E_{|U_i}\to V_i \times\R,$
where $V_i$ denotes a bounded open subset of $\R^n$, such that
\begin{enumerate}
\item $\forall i\in I$, there exist smooth kernels
$k_i \in 
\Gamma(M\times M, E^*\boxtimes E) $
such that for every $s_i\in \Gamma(M,E)$ with support in $U_i$
and every $x\in M\setminus U_i$, 
$$\ P(s_i)(x) = \int_M k_i(x,y) s_i(y) |dy|,$$
where $|dy|$ denotes a Lebesgue measure on $M$.
\item $\forall i\in I,$  there exist smooth symbols $p_i : V_i \times \R^n \cong
T^*M_{|U_i} \to \C$ such that for every $s_i \in \Gamma(M,E)$
with support in $U_i$ and every $x\in V_i$,
$$ \Phi_i(P(s_i))(x)= \iint_{V_i \times \R^n } p_i(x,\xi)e^{i\cg x-y,\xi\cd} \Phi_i(s_i)(y)  d\xi dy,$$
where $dy d\xi$ si the standard Lebesgue measure on $V_i \times \R^n $.
\item For every compact subset $K_i \subset V_i$ and every $\alpha, \beta \in \Nn^n$,
there exist positive constants $c_{K_i, \alpha, \beta}$ such that 
$$ \forall (x,\xi)\in K_i\times \R^n , 
|\frac{ \partial}{\partial x^\beta} \frac{\partial }{\partial \xi ^\alpha} p_i (x,\xi)| \leq c_{K_i, \alpha, \beta} (1+|\xi|)^{m-|\alpha|}.$$
\end{enumerate}
\end{defin}
Now, let $h_E$ be a Riemannian metric on $E$ and $|dy|$
be  a Lebesgue measure on $M$, 
which we assume to be  compact and without boundary. Then, $\Gamma (M,E)$
inherits the $L^2$-scalar product (\ref{<>}).
\begin{defin}\label{adj}
The adjoint of the pseudo-differential operator $P$ is the operator ${}^t P$
satisfying for every 
$ s,t\in\gme, \cg P(s),t\cd = \cg s, {}^t P(t)\cd.$
When ${}^t P = P,$ the operator is said to be self-adjoint.
\end{defin}

\begin{defin}(see \cite{Hormander-psd}, \cite{Hormander}, \cite{Hormander-Fourier})\label{def 2.2}
A self-adjoint pseudo-differential operator of order $m\in \R$ given by Definitions \ref{Definition 4.0},
\ref{adj}
is said to be elliptic if and only if for every $i\in I$ and 
every $(x,\xi)\in T^*M_{|U_i}$ such that $\xi \not=0$, 
the limit
 $$\si_P(x,\xi) = \lim_{t\to + \infty} \frac{1}{t^m} p_i(x,t\xi)$$
exists and is positive. This limit then does not depend on the choice
of $i\in I$ and defines a positive homogeneous  function  
$\si_p : T^*M \to \R$ of order $m$ and class $C^\infty$. 
\end{defin}
The function $\sigma_P$ given by Definition \ref{def 2.2}
will be called the homogenized 
principal symbol of $P$. It is symmetric in the sense that 
for every 
$(x,\xi)\in T^*M, \ \sigma_P(x,-\xi) = \sigma_P(x,\xi).$

\begin{exa}
Recall that if in a local trivialization of $E$ the differential 
operator $Q$ of order $m$ 
reads 
$
f\in C_c^\infty(\R^n , \R) \mapsto \tilde Q (\partial /\partial x_1, \cdots, \partial /\partial x_n)(f)\in C^\infty_c(\R^n, \R),$
where $\tilde Q \in C^\infty(\R^n) [X_1, \cdots, X_n]$, 
and if $\tilde Q_m$ is the homogeneous part of order $m$ of $\tilde Q$, 
then the principal symbol of $Q$ 
is the homogeneous function of order $m$
$
\sigma_Q : (\R^n)^*\to \C $ satisfying
$\si_Q (\xi_1 dx_1 + \cdots + \xi_n dx_n) = \tilde Q_m (i\xi_1, \cdots, i\xi_n).
$
\end{exa}

\begin{defin}
An elliptic self-adjoint pseudo-differential operator $P$ on $\gme$ is said 
to be bounded from below if and only if there exists a constant 
$c\in\R$ such that for every 
$  s\in\gme, \cg P(s),s\cd \geq c \cg s,s\cd.$
It is said to be positive when $c > 0$.
\end{defin}
\begin{rem}\label{positif}
The transformation $P \to P- c Id $ turns any elliptic self-adjoint pseudo-differential operator
bounded from below into a positive one. Since our results are not sensitive to this transformation,
they hold for  any operator bounded from below even if we sometimes assume it
to be positive for simplicity.
Recall finally that these operators have discrete spectrum with finite dimensional eigenspaces. 
\end{rem}

\subsection{Proof of Theorem  \ref{Theoreme 3} }\label{A1}

Set $L= \lambda^m$ and $\el= e_L$. 
The strategy followed by H\"ormander is the following.
The derivative of $\el$  with respect to $\lambda$
is a distribution whose support 
is the set of eigenvalues of $P$. 
Its Fourier transform with respect to $\lambda$
is the kernel of the hyperbolic equation 
$\partial_t u + i P^{1/m}= 0$, where $P^{1/m}$
stands for the operator with the same eigenfunctions as $P$
and whose eigenvalues are the $m$-th root
of the corresponding ones of $P$. H\"ormander proves that 
in a neighbourhood $V$ of the diagonal of $M\times M$
and for small values of the time $t$,
this kernel takes the form of a Fourier integral operator, modulo an operator with 
 smooth
kernel.
Consequently, if $\rho : \R \to \R$ 
is a non negative function in the Schwartz space
such that 
its Fourier transform $\hat \rho$ satisfies
$\hat \rho(0)= 1$ and $Supp(\hat \rho) \subset
[-\epsilon, \epsilon]$, 
then for every $x,y\in V$,
$$ 
\int_{-\infty}^{+\infty} \rho(\lambda-\mu)
\partial_{\mu} \tilde e_\mu (x,y) d\mu 
- \int_{T_y^*M} R(x,\lambda - p'_{|y}(\xi'), y, \xi)
e^{i\psi(x,y,\xi)} d\xi$$
is a  rapidly decreasing  function as $\lambda \to +\infty$,
where 
\begin{itemize}
\item $\psi(x,y,\xi)= \cg x-y, \xi \cd
+ O(|x-y|^2|\xi|)$ when $x\to y$,
for a scalar product $\cg \, , \, \cd$ 
in a chart of $M$ that contains $x$ and $y$.
\item $p'(\xi)= \sigma_P (\xi)^{1/m} + O(1)$
\item $R(x,\lambda,y,\xi)=\frac{1}{2\pi}
\int_\R \hat \rho(t) q(x,t,y,\xi) 
e^{it\lambda} dt$
with 
$q(x,0,y,\xi)= (\frac{1}{2\pi})^n + O(1/|\xi|),$
see Lemma 4.1 of \cite{Hormander}.
\end{itemize}
This function $R$ is rapidly decreasing as $\lambda$
grows to infinity.
After differentiation we deduce likewise that
\bq \label{Q1Q2}
\int_{-\infty}^{+\infty} \rho(\lambda-\mu)
\partial_{\mu} Q_1Q_2\tilde e_\mu (x,y) d\mu 
- \int_{T_y^*M} Q_1Q_2 (R(x,\lambda - p'_{|y}(\xi), y, \xi)
e^{i\psi(x,y,\xi)} )d\xi
\eq
is a rapidly decreasing function as $\lambda $ 
grows to infinity.
\begin{lem}\label{Lemme 7}(Compare
Lemma 4.3 of \cite{Hormander})
Under the hypotheses of Theorem \ref{Theoreme 3},
there exists a constant $c>0$
such that for every $(x,y)$ in a neighbourhood
$V$ of the diagonal of $M\times M$, for every 
$\lambda \geq 0$ and every $0\leq \mu \leq 1$,
$$ 
\| Q_1Q_2 \tilde e_{\lambda +\mu}(x,y)
- Q_1Q_2 \tilde e_{\lambda }(x,y)\|_{h_E}
\leq
C(1+|\lambda|)^{n-1+|\sigma_{Q_1}|+ |\sigma_{Q_2}|}.
$$
\end{lem}
\bpr
Let us assume first that $Q_1=Q_2$ and $x=y$.
We proceed as in the proof of Lemma 4.3 
of \cite{Hormander}.
The function 
$$\partial_\mu Q_1 Q_1 \tilde e_\mu (x,x)
= 
\sum_{k} \delta_{\lambda_k}
h_E(Q_1s_k(x), Q_1s_k(x))$$
is positive, where $s_k$ 
is an eigenfunction with eigenvalue $\lambda_k^m$.
We deduce the existence of a constant 
$C_1>0$
such that 
$$\|  Q_1Q_1 \tilde e_{\lambda +\mu}(x,y)
- Q_1Q_1 \tilde e_{\lambda }(x,y)\|_{h_E} 
\leq
C_1\int_\R  \rho(\lambda-\mu)\partial_\mu Q_1Q_1 \tilde e_\mu (x,x) d\mu.$$
From (\ref{Q1Q2}), it is enough to bound from above the integral
$$ \int_{T_x^*M} Q_1Q_1 \big(R(x, \lambda - p'(\xi), y,\xi )e^{i\psi(x,y,\xi)}\big)d\xi$$
From the ellipticity of $P$ we deduce the existence of $C_2>0$
such that 
$$ \forall \xi \in T_x^*M, | Q_1Q_1 \psi(x,y,\xi)|_{|(x,x)}=
|\sigma_{Q_1}(\xi)|^2 \leq C_2 (1+ p'(\xi))^{2|\sigma_{Q_1}|}.$$
Following \cite[p 210]{Hormander}, we deduce that 
\beq
| \int_{T_x^*M} Q_1Q_{1|(x,x)} (R(x,\lambda - p'(\xi), y, \xi)) e^{i\psi(x,y,\xi)}d\xi |
& \leq & 
C_3\int_\R (1+ |\lambda - \sigma|^{-N})(1+ |\sigma|)^{2|\sigma_{Q_1}|} dm(y,\sigma)\\
& \leq & O(\lambda^{-\infty}) + C_4 (1+ |\lambda|)^{n-1+2|\sigma_{Q_1}|},
\eeq
where $C_3$, $C_4$ are positive constants, $N$ denotes 
a large enough integer  and 
$$m(x,\sigma)= \int_{\{\xi \in T_x^*M| \ \sigma_P(\xi)\leq \sigma\}}
d\xi.$$
We deduce the result when $Q_1=Q_2$ and $x=y$, then likewise
when $(x,y)$ lies in a neighbourhood $V$
of the diagonal, see Lemma 3.1 of \cite{Hormander}. 
The general case is now a consequence of the Cauchy-Schwarz inequality
and there exists a positive constant $c$ such that 
$ \forall x,y\in N, \forall \lambda >0, \forall \mu \in [0,1],$
\beq
\| Q_1Q_2 \tilde e_{\lambda +\mu}(x,y)
- Q_1Q_2 \tilde e_{\lambda }(x,y)\|_{h_E} & = &
\| \sum_{k \ | \ \lambda \leq \lambda_k \leq \lambda + \mu}
Q_1(s_k(x))Q_2(s_k(y))\| \\
&  \leq & 
\big(\sum_{k \ | \ \lambda \leq \lambda_k \leq \lambda + \mu}
\| Q_1(s_k(x))\|^2\big)^{1/2}\cdots\\
&&\cdots \big( \sum_{k \ | \ \lambda \leq \lambda_k \leq \lambda + \mu}
\| Q_2(s_k(y))\|^2\big)^{1/2} \\
& \leq &( 
\| Q_1Q_{1|(x,x)}\tilde e_{\lambda  +\mu} - 
Q_1Q_{1	|(x,x)}\tilde e_{\lambda }\|^2)^{1/2}\cdots\\
& &\cdots
( 
 \| Q_2Q_{2|(y,y)}\tilde e_{\lambda  +\mu} - 
Q_2Q_{2|(y,y)}\tilde e_{\lambda }\|^2)^{1/2}\\
& \leq & C(1+|\lambda|)^{n-1+|\sigma_{Q_1}|+|\sigma_{Q_2}|}.
\eeq
\epr
\bpr[ of Theorem \ref{Theoreme 3}]
We proceed as in \cite{Hormander}, p. 211. 
We deduce from Lemma \ref{Lemme 7}
that $ \forall x,y\in U, \forall \lambda \geq 0, \forall \mu \geq 0, $
$$ 
\| Q_1Q_2 \tilde e_{\lambda + \mu}(x,y)
- Q_1Q_2 \tilde e_\lambda (x,y)\|_{h_E}\leq
C(1+\lambda + \mu)^{ n-1+ |\sigma_{Q_1}| + |\sigma_{Q_2}| }(1+\mu).$$
Thus, there exists $C'>0$ such that
$$ \| \int_\R \rho(\lambda-\mu)Q_1Q_2 \tilde e_\mu(x,y)d\mu-
Q_1Q_2 \tilde e_\lambda (x,y)\|_{h_E} \leq
C' (1+ \lambda)^{n-1+ |\sigma_{Q_1}|+|\sigma_{Q_2}|}.$$
However, by integration of (\ref{Lemme 7})
over the interval $]-\infty, \lambda]$, 
we deduce the existence of $C''>0$ such that
$$ \| Q_1Q_2 \tilde e_{\lambda + \mu}(x,y)
- \int_{T_y^*M} \int_{-\infty}^\lambda
Q_1 Q_2(R(x,\sigma - p'_{|y}, y, \xi) e^{i\psi (x,y,\xi)})
d\xi d\sigma\|_{h_E} \leq C''.$$
Moreover, by definition of $\psi$ and $R$,
$\int_{T_y^*M} \int_{-\infty}^\lambda
Q_1 Q_2(R(x,\sigma - p'_{|y}(\xi), y, \xi) e^{i\psi (x,y,\xi)})d\xi d\sigma$ equals
\beq
\frac{1}{(2\pi)^n} \int_{\{\xi \in T_y^*M\, | \, p'(\xi )\leq \lambda\}}
(1+ O(1/|\xi|))Q_1Q_2 e^{i\psi(x,y,\xi)} d\xi 
+\cdots \\
\cdots \int_{T_y^*M} Q_1 Q_2 (R_1 (x,\lambda - p'_{|y}(\xi), y , \xi) e^{i\psi(x,y,\xi)} ) d\xi,
\eeq
where
$$
R_1 =  \left\{ \begin{array}{cc} \int_{-\infty}^\tau R(x,\si,y,\xi) d\sigma & \text{ if } \tau \leq 0 \\
\int_{-\infty}^\tau R(x,\si,y,\xi) d\sigma - q(x,0,y,\xi)
& \text{ if } \tau >0
\end{array}\right.
$$
is a function which decreases faster than any polynomial,
see \cite[p. 211]{Hormander}. 
Thus, there exists a constant $C'''>0$ such that 
\beq
\| \int_{T_y^*M \times ]-\infty, \lambda] }
Q_1 Q_2 (R(x,\sigma - p'_{|y} (\xi), y, \xi ) e^{i\psi (x,y,\xi)})d\xi d\sigma - \cdots & &\\
\cdots \frac{1}{(2\pi)^n} \int_{\xi \in T_y^*M \, | \, p'(\xi )\leq \lambda \} } \sigma_{Q_1} (\xi)\overline{\sigma_{Q_2}(\xi)} d\xi\| 
&\leq & C'''(1+\lambda )^{n-1+ |\sigma_{Q_1}| + |\sigma_{Q_2}| }.
\eeq
From the triangle inequality, we finally deduce 
that there exists $C''''>0$ such that
for every $(x,y)\in V$, 
$$ \| Q_1 Q_2 e_L(x,y)- 
\frac{1}{(2\pi)^n} \int_{\{\xi \in T_y^*M \, | \ \sigma_P (\xi) \leq L\}}
\sigma_{Q_1}(\xi)\overline{\sigma_{Q_2}(\xi)} d\xi \| \leq 
C''''(1+ \lambda)^{n-1+ |\sigma_{Q_1}| + |\sigma_{Q_2}|}.
$$
Hence the result.
\epr

%\end{enumerate}
%\end{rem}
% % % % % % % % % % % % % % % % % % % % % % % % % % % % % % % % %
\normalsize
\baselineskip=17pt
%%%%%%%%%%%%%

\providecommand{\bysame}{\leavevmode\hbox to3em{\hrulefill}\thinspace}
\providecommand{\MR}{\relax\ifhmode\unskip\space\fi MR }
\providecommand{\MRhref}[2]{%
  \href{http://www.ams.org/mathscinet-getitem?mr=#1}{#2}
}
\providecommand{\href}[2]{#2}

Damien Gayet\\
 Univ. Grenoble Alpes, IF, F-38000 Grenoble, France\\
CNRS, IF, F-38000 Grenoble, France\\
 damien.gayet@ujf-grenoble.fr\\

Jean-Yves Welschinger\\
 Universit\'e de Lyon \\
CNRS UMR 5208 \\
Universit\'e Lyon 1 \\
Institut Camille Jordan \\
43 blvd. du 11 novembre 1918 \\
F-69622 Villeurbanne cedex\\
France\\
 welschinger@math.univ-lyon1.fr

\end{document}